\documentclass[11pt,a4paper,amssymb,amsmath, tightenlines]{article}

\usepackage{amsmath, amssymb, amsthm, latexsym, mathrsfs, url, color, epsfig, graphics, float, setspace,hyperref, graphicx,bbm}
\usepackage{tikz}
\usepackage{sectsty}
\usepackage{float}
\usepackage[T1]{fontenc}
\usepackage[bitstream-charter]{mathdesign}
\usepackage{pgfplots}
\usepackage{color}
\usepackage{bm}
\usepackage{soul}
\usepackage[numbers]{natbib}
\usepackage{algorithmicx}
\usepackage{algorithm}
\usepackage{algpseudocode}
\usepackage{appendix}

\usepackage{enumerate}
\usepackage{tikz}
\usepackage[utf8]{inputenc}
\usepackage{pgfplots} 
\usepackage{pgfgantt}
\usepackage{pdflscape}
\newlength\fwidth
\newlength\fheight

\usepackage{multirow}

\usepackage{caption}
\usepackage{subcaption}

\pgfplotsset{compat=newest} 
\pgfplotsset{plot coordinates/math parser=false}

\allowdisplaybreaks[1]

\newtheorem{Theorem}{Theorem}[section]
\newtheorem{Definition}{Definition}[section]
\newtheorem{rem}{Remark}
\newtheorem{Proposition}{Proposition}[section]

\newtheorem{Corollary}{Corollary}[section]
\newtheorem{Example}{Example}[section]
\numberwithin{equation}{section}

\newcommand{\lnn}{{\rm ln}}
\newcommand{\expp}{{\rm exp}}
\newcommand{\erf}{{\rm erf}}

\newcommand{\F}{\mathscr{F}}

\newcommand{\PR}{\mathbb{P}}

\newcommand{\Q}{\mathbb{Q}}
\newcommand{\R}{\mathbb{R}}

\newcommand{\rd}{\textup{d}}

\newcommand{\indi}[1]{1\hspace{-.09cm}\textup{\textrm{l}}}

\begin{document}
\title{\vspace{-2cm}\bf Quantile diffusions for risk analysis}
\author{Holly Brannelly$^{1}$, Andrea Macrina{$^{1,* }$\footnote{$*$ Corresponding author: a.macrina@ucl.ac.uk}}, Gareth W. Peters{$^{2}$} \\ \\ {$^{1}$Department of Mathematics, University College London} \\ {London, United Kingdom} 
\\ {$^{2}$}Department of Statistics \& Applied Probability \\ University of California Santa Barbara \\ Santa Barbara, USA}
\date{10 September 2021}
\maketitle
\vspace{-1cm}
\begin{abstract}
\noindent 
\\\vspace{-1cm}\\
We develop a novel approach for the construction of quantile processes governing the stochastic dynamics of quantiles in continuous time. Two classes of quantile diffusions are identified: the first, which we largely focus on, features a dynamic random quantile level and allows for direct interpretation of the resulting quantile process characteristics such as location, scale, skewness and kurtosis, in terms of the model parameters. The second type are function--valued quantile diffusions and are driven by stochastic parameter processes, which determine the entire quantile function at each point in time. By the proposed innovative and simple---yet powerful---construction method, quantile processes are obtained by transforming the marginals of a diffusion process under a composite map consisting of a distribution and a quantile function. Such maps, analogous to rank transmutation maps, produce the marginals of the resulting quantile process.  We discuss the relationship and differences between our approach and existing methods and characterisations of quantile processes in discrete and continuous time. As an example of an application of quantile diffusions, we show how probability measure distortions, a form of dynamic tilting, can be induced. Though particularly useful in financial mathematics and actuarial science, examples of which are given in this work, measure distortions feature prominently across multiple research areas.  For instance, dynamic distributional approximations (statistics), non--parametric and asymptotic analysis (mathematical statistics), dynamic risk measures (econometrics), behavioural economics,  decision making (operations research), signal processing (information theory), and not least in general risk theory including applications thereof, for example in the context of climate change.
\\\vspace{-0.5cm}
\end{abstract}

\section{Introduction}
The dynamics of a stochastic process can be characterised in numerous ways, including through trend, volatility, higher--order moment dynamics, and under transformations of the finite dimensional distributions or transition law of the process.  In this work we develop quantile diffusions, a novel class of stochastic processes in continuous time. We then identify the class of Tukey quantile diffusions of which the parametrisation enables direct interpretation with regard to higher--order moments such as skewness and kurtosis. The quantile processes developed in this paper accommodate a wide range of tail behaviours in the finite dimensional distributions, ranging from exponential decay to heavy--tailed regular variation.

The concept of a measure distortion framework, often applied on the density space, is ubiquitous in statistical science under the topics of measure and density approximation. In such settings density distortions are used to expand, modify or tilt a base distribution such that the distortion alters the moments or cummulants of the distribution relative to the base distribution in such a fashion that the resulting distorted distribution may better satisfy a target objective, which is often expressed in terms of moments or cummulants, thereby producing an improved approximation. Common examples include the saddle point, Edgeworth and generalised Esscher transformations. Such approaches have enjoyed widespread use to produce transformations applied to a family of base densities in order to distort the moment or tail behaviour characteristics of the resulting distribution, relative to the base distribution, for some analytic purpose. For an overview of this large family of methods see, e.g.,  \cite{barndorff1979edgeworth,bickel1974edgeworth,wand1991transformations,daniels1954saddlepoint}. In this paper we wish to characterise a family of distortions, which extend the static case to a dynamic, continuous--time setting. Furthermore, rather than working with the density space, we wish to achieve the dynamic distortions directly in the quantile space, since many risk management problems are expressed in terms of quantile functions or quantile processes rather than a density. As such, one way to express the objective of this paper is the construction of a very flexible family of transformations, with which---when applied to a Markov stochastic process---one can produce a dynamic, continuous--time family of distortion processes induced by a family of quantile processes, characterised by the parametric form of the selected transformation composite map and the underlying base diffusion process. The utility of quantile dynamics in this context is under--explored and to our knowledge it has not yet been proposed or studied in the setting of stochastic diffusions.

The advantage of working with distorted measure flows in the quantile space, as opposed to the density space, is that many applications in risk management require the quantification of the risk according to quantile functions, rather than the distorted density. In this work we are able to demonstrate how to consistently and rigorously transform a process to a quantile process in continuous time that serves the purpose of producing a measure flow of distorted quantiles, which will consequently be interpretable in its properties. It is the purpose of this paper to lay the mathematical foundations for this framework to be rigorously defined, and then to demonstrate through applications how general such a framework could be for applications in a variety of disciplines with focus on risk analysis including, but by no means limited to, financial mathematics and actuarial science. 

We emphasize that such a methodology can be adopted for a variety of applications including dynamic risk measures in econometrics, e.g., the study of VaR processes, behavioural economics and dynamic consumer preference models, operations research and dynamic decision making, signal processing in information theory, as well as to more general applications requiring dynamic risk measures such as disaster monitoring (earthquake hazard, flood hazard and other natural disasters). We also draw attention to dynamic risk analysis in climate and environmental science where the risk measure is not a monetary figure, necessarily. As such the general mathematical framework proposed in this manuscript will readily translate into such applications and many other.

Dynamical quantile functions, quantile processes and quantile diffusions have been explored previously in the stochastic process literature, the statistical regression and time series econometrics literature, the risk management and insurance literature, as well as the mathematical statistical literature within the study of empirical processes. As a result, there are numerous meanings attributed to the terminology ``quantile process'' or ``quantile dynamics'' based on the definitions developed in earlier works. In Section \ref{qdconstruction} we briefly identify and compare these existing characterisations in order to illustrate the novelty of our proposed framework. A comprehensive discussion on the use of quantile functions in discrete, time series based statistical modelling and data analysis is given in \cite{gilchrist, koenkerbassett} and a tutorial review in \cite{quantiletutorial}. In the risk management setting, the literature on dynamical risk measures parallels ideas presented in the statistical and time series econometrics literature when it comes to dynamical risk measure processes---see examples in both discrete \cite{bielecki2017survey,  chavez2014extreme,dqm,quantiletutorial,lmomentsgh} and continuous time contexts \cite{acciaio2011dynamic,bion2008dynamic,detlefsen2005conditional, riedel2004dynamic}. Distinct from these statistical time series modelling and dynamic risk management frameworks, there have also been developments of what are termed ``quantile processes'' for empirical processes in mathematical statistics and probability literature---see for instance the sequence of works \cite{qp1,qp2, qp3}. 

Since the focus of the work in this paper pertains to modelling in continuous time, we focus on the literature on diffusion processes where quantile processes have also been explored.  Here, one may consider the works of \cite{akahorioption, dassios, embrechts1995proof} and \cite{yor} whereby---building on ideas in \cite{miura}---a Brownian process is considered and, at each instance in time, the distribution of a random variable, defined by the $\alpha$--quantile of the diffusion at this time, is studied. 

From a financial mathematics perspective, Miura \cite{miura} motivated such a consideration by introducing the ``$\alpha$--percentile option'' whereby the underlying is given by the $\alpha$--percentile of the price process over the life of the option, e.g., the median if $\alpha=0.5$.  Knowing the distribution of the $\alpha$--quantile allows for the pricing of these path--dependent options, as shown in \cite{akahorioption}.  Whilst the models constructed in this paper may also be used to price options written on a quantile, comparatively to work where the focus is the distributional behaviour of the ``quantiles of diffusions'' (the quantile dynamics of the diffusion process is studied), our work focuses on producing a new class of quantile diffusion processes. Whilst an underlying diffusion model is implied, rather than specifying this and then working on characterising the quantiles as is classically done in the literature, we work directly on constructing a process in the space of quantile functions.  This can be scalar--valued or function--valued depending on which form of our approach one utilises. There are several advantages to this alternative approach that we explain in the remainder of the paper.

We develop two classes of quantile diffusions, one based on a ``process--driven'' (or ``random--level'') construction, which produces a scalar--valued process, and the second based on a ``parameter--driven'' (or ``function--valued'') construction, which produces a function--valued process in the space of quantile functions. In this paper we largely focus on the first of these constructions, where we develop a family of quantile diffusion processes by transforming each marginal of a given univariate diffusion process under a composite map consisting of a distribution function and a quantile function, which in turn produces the marginals of the resulting quantile process. This is analogous to the rank transmutation mapping approach, see \cite{distalchemy, rtm, quantilemech2, quantilemech}. We utilise a flexible family of such maps that allows for the moments of the underlying process to be directly interpreted with regard to the parameters of the transformation, based on the Tukey $g$--and--$h$ transformation map, see \cite{tukey}. In general, the notion of a transmutation map is to transform from some ``base'' distribution $F$ to a ``target'' distribution with quantile function $Q$, as discussed in \cite{quantilemech2}. Here, a differential equation that is referred to as the ``recycling ODE'' is derived, of which solutions provide a direct route to the object $G(x)=Q(F(x))$ when the inverse of the distribution functions $F$ and $Q$ may not be easily available. The only requirement here is the ability to calculate the logarithmic derivatives of the two corresponding densities. The motivations behind such a map involve the ability to provide a one--step approach to introduce features such as skew or kurtosis into a distribution that may, for example, be symmetric, in order to more realistically model financial asset returns. Similarly, sample transmutation maps and rank transmutation maps are explored in \cite{distalchemy, rtm}, and in \cite{gilchrist} (where these are referred to as $Q$--transformations and $P$--transformations, respectively), again providing a succinct method of moving from the distributional setting to the quantile setting whilst introducing relative skewness or kurtosis into a given distribution without the use of Gram--Charlier expansions (which can be viewed as the asymptotic analogue of the rank transmutation maps).  These maps also allow for converting samples from one distribution into those from another without the need for Cornish--Fisher expansions. A quadratic rank transmutation map is used in \cite{rtm} to produce skew--uniform, skew--exponential and skew--normal distributional representations. A brief comparison to existing literature revolving around the idea of modulating a distribution to introduce skewness, such as that by Azzalini, see \cite{azzalini2, azzalini1,azzalini3}, and Genton, see \cite{gentonskew}, is made.  The advantages of distortion maps that produce quantiles over models of the Azzalini type lies in the ability to introduce relative skewness to some base distribution, as opposed to an absolute amount of skewness, thus providing substantial practical flexibility since any base--line model may be used.  

In this paper the class of maps is selected to directly determine the stochastic dynamics (which satisfy an SDE) of the quantiles, or entire quantile curves, through time.  This construction--based approach allows us to obtain ``target'' quantile models, with a high level of flexibility over the statistical properties of the resulting model. The emphasise in our paper is on the {\it construction approach} we propose, by which new classes of quantile diffusion processes can be constructed explicitly. This is in contrast to, say, an approach that defines and characterises continuous--time quantile diffusions by studying their general properties but stops short of providing a ``recipe'' for building such processes.

\subsection{Contributions}
The aim of this work is to extend the discrete--time quantile models and quantile transformation maps to the continuous--time diffusion setting.  We develop two approaches of constructing such models, the first of which distorts each marginal of some given process through a composite map that is of similar form to a rank transmutation map.  Quantile diffusions are generated by applying this mapping to a scalar--valued stochastic diffusion (the ``driving process''), the output being a scalar--valued quantile diffusion. We emphasise the wide class of models with directly interpretable statistical characteristics that arise from such a transformation. The resulting quantile diffusion process satisfies an SDE for the dynamics of the quantiles. The properties of the model, captured in the drift and volatility functions of the process, then depend entirely on the choice of functions involved in the mapping, the drift and volatility coefficients, the autocorrelation and other types of extreme tail dependence structures of the driving diffusion.

Each realisation of the output process corresponds to a single quantile level, and hence when the paths of the underlying driving process are sampled ``infinitely'' many times, and the samples at each fixed time are ordered, samples of the quantile diffusion representing all quantile levels in $[0,1]$ at that time are obtained. Producing these ordered samples at each time $0< t<\infty$ allows one to model the time evolution of the entire quantile function. Once the new class of quantile diffusion processes has been defined, one has the same level of flexibility as with the quantile--preserving maps discussed in \cite{quantiletutorial}, and hence the properties of the model can be chosen so as to alter the symmetry or tail properties of the process. 

The second approach is a parameter--driven model whereby we put a multivariate diffusion on the parameters of a well--defined quantile function, and hence map from realisations of these parameter processes to function--valued realisations of the quantile diffusion.  Each sample path of the (possibly multidimensional) parameter process will drive the resulting function--valued quantile diffusion, allowing one to dynamically model the entire quantile function at any instance in time.  A discrete--time equivalent to this model is given in \cite{dqm}. 

In Section \ref{qdconstruction} we discuss different approaches to quantile diffusion constructions, contrasting the differences between such approaches and our approach. In Section \ref{QuantDiffSDE} we derive the SDEs, detailing specific conditions for existence and uniqueness of the solutions and characterising the strong and weak solutions. In Section \ref{familiesofmodelssection} we apply this approach to a variety of sub--families of Tukey models.  In Section \ref{examplessection} we consider an application of our approach in the context of ``distortion--based pricing'' found in financial and insurance mathematics. Some proofs and additional technical material are included in the appendix.  

\subsection{Notation and definitions}
The notation and common definitions used in this paper includes the formulation of the relevant probability spaces. Let $(\Omega,\mathscr{F}, (\mathscr{F}_t)_{0\leq t<\infty},\mathbb{P})$ denote a filtered probability space with filtration $(\mathscr{F}_t)$ and $(W_t)_{0\leq t<\infty}$ an $(\mathscr{F}_t)$--adapted, one--dimensional, standard Brownian motion. 

We will define the space of quantile functions for a random variable $X$ with distribution $F_X$ as given in Definition \ref{quantilefunctiondef}, which utilises the generalised inverse discussed in \cite{geninverse} and provided in Definition \ref{generalisedinversedef}. This generalised inverse allows one to consider instances where one wants to invert a distribution function that may not be real--valued, continuous and strictly monotone, and hence the ordinary definition of the inverse it possesses on its range does not apply.

\begin{Definition}\label{generalisedinversedef}
For an increasing function $F:\mathbb{R}\rightarrow\mathbb{R}$ with $F(-\infty)=\lim_{x\downarrow -\infty}F(x)$ and $F(\infty)=\lim_{x\uparrow\infty}F(x)$, the collection $\Q$ of generalised inverse functions $Q:=F^{-}:\mathbb{R}\rightarrow[-\infty,\infty]$ of $F$ is defined by
\begin{equation}\label{geninverse}
Q(y)=\inf\left\{x\in\mathbb{R}:F(x)\geq y \right\}, \hspace{2mm} y\in\mathbb{R},
\end{equation}
with the convention that $\inf\emptyset=\infty$. 
\end{Definition}  
In what follows, we denote the distribution, quantile and density functions of a random variable by $F_{\bullet},Q_{\bullet},$ and  $f_{\bullet}$, respectively, where the argument in the subscript of such functions denotes the random variable to which they correspond. We denote the collection of quantile functions of the random variable by $\mathbb{Q}_{\bullet}$ with the same argument in the subscript.  When included in notation, parameters will follow a semi--colon in the arguments of these functions.
\begin{Definition}\label{quantilefunctiondef}
Let $X$ be a real--valued random variable with distribution function $F_X:\mathbb{R}\rightarrow[0,1]$.  The corresponding quantile function of $X$ is $Q_X=F_X^{-}:[0,1]\rightarrow[-\infty,\infty]$ where $Q_X\in\mathbb{Q}_X\subset\mathbb{Q}$.
\end{Definition}

\begin{Definition}\label{mappingclassdefn}
Consider a space of functions $\mathbb{S}\subseteq\mathbb{Q}_{X}\subset\mathbb{Q}$. We define the class of functions that characterise such a space by $m:\mathbb{Q}\rightarrow\mathbb{S}$ for $m$ a continuous, increasing, bijective mapping.
\end{Definition} 

We introduce a generic diffusion process in continuous time as follows:
\begin{Definition}\label{diffusiondef}
A diffusion is a process $(Y_t)_{0\leq t<\infty}$ on the probability space $(\Omega,\mathscr{F}, (\mathscr{F}_t)_{0\leq t<\infty},\mathbb{P})$ satisfying
\begin{equation}\label{ysde}
\rd Y_t =  \mu\left(t,Y_t\right)\rd t + \sigma\left(t,Y_t\right)\rd W_t
\end{equation}
where $Y_0=y_0\in\mathbb{R}$, $\mu(t,y):\mathbb{R}^+\times\mathbb{R}\rightarrow\mathbb{R}$ and $\sigma(t,y):\mathbb{R}^+\times\mathbb{R}\rightarrow\mathbb{R}^+$ satisfy the necessary conditions, see e.g., \cite{oksendal}, to ensure the SDE (\ref{ysde}) admits a solution $Y_t = Y(t,\omega):\mathbb{R}^+\times\Omega\rightarrow\mathbb{R}$.    Let the transition distribution of $(Y_t)$ be given by $F_Y\left(t,y \vert\, s,x\right):=\mathbb{P}(Y_t\leq y \vert\, Y_s=x )$ for $0\leq s\leq t\leq T<\infty$ and $x, y\in\mathbb{R}$. The transition density function is given by $f_Y\left(t,y\vert\,s,x \right) := \partial_y F_Y\left(t,y \vert\, s,x \right)$ if the derivative exists.  
\end{Definition}
A statement on the existence of weak and strong solutions of SDEs is given in Appendix \ref{strongsolutionsection}.
\section{Construction of quantile diffusions}\label{qdconstruction}
In this section we introduce our approach to constructing stochastic quantile processes in continuous time. We highlight differences and distinguish between our approach proposed in this paper and existing approaches to studying quantile dynamics.
\subsection{Characterisations of quantile processes}\label{quantileprocesses}
We recall three approaches for the definition of quantile processes, in order to differentiate them from the use of this terminology in our constructions. In much of the literature the formulation adopted in \cite{qp2} is invoked when one refers to a quantile process, and it is based on the univariate, empirical quantile process given in Definition \ref{quantileprocessdefn}.

\begin{Definition}\label{quantileprocessdefn}
Let $Y_1,Y_2,\ldots, Y_n$ be a sequence of i.i.d. random variables with a continuous distribution function $F_Y$, and let $Y_{(1,n)}\leq Y_{(2,n)}\leq \ldots \leq Y_{(n,n)}$ denote the order statistics of the random sample $Y_1,Y_2,\ldots,Y_n$.  Define the empirical distribution function $F_n(y)$ and the quantile function $Q_n(u)$ as follows:
\begin{align*}
F_n(y) &= \begin{cases}
0 &{\rm{if}} \hspace{2mm} Y_{(1,n)}>y \\
\frac{k}{n} &{\rm{if}} \hspace{2mm} Y_{(k,n)}\leq y<Y_{(k+1,m)}, \hspace{2mm} k=1,2,\ldots,n-1 \\
1 &{\rm{if}} \hspace{2mm}Y_{(n,n)}\leq y
\end{cases}
\\[10pt]
Q_n(u) &=  Y_{(k,n)} \hspace{2mm} {\rm{if}} \hspace{2mm} \frac{k-1}{n}< u\leq \frac{k}{n}, \hspace{2mm} k=1,2,\ldots,n.
\end{align*}

The empirical quantile process is defined by $q_n(u) = n^{1/2}\left(Q_n(u)-F_Y^{-}(u)\right)$ for $0<u<1$.
\end{Definition}

This definition relates to the convergence of the law of the order statistics of an empirical process, which is observed as a sequence of independent and identically distributed (i.i.d.) random variables from a fixed distribution $F_Y$. 

The second widely adopted definition of a quantile process, this time in the context of a quantile diffusion (in particular the quantile of a Brownian motion with drift) is given as follows, see \cite{dassios, embrechts1995proof}.  The extension to processes with stationary and independent increments is given in \cite{dassiosquantile2}. 

\begin{Definition}\label{quantileprocessdefn2}
Let $(\Omega,\mathscr{F},(\mathscr{F}_t),\mathbb{P}^{\mu})$ be a probability space where $\mathbb{P}^{\mu}$ denotes the law of the Brownian motion with drift $(Y_t)_{0\leq t<\infty}$, given by $Y_t:=W_t+\mu t$ for $\mu\in\mathbb{R}$, on the canonical space $(C(\mathbb{R}_+, \mathbb{R}),\mathscr{F}_{\infty})$, and where $\mathscr{F}_t=\sigma((Y_s)_{0\le s\le t})$, $t\in[0,\infty)$.  For $\alpha \in [0,1]$ and $\omega\in\Omega$, define the $\alpha$-quantile diffusion $(M_t^{(\alpha)})_{t>0}$ by 
\begin{equation}
M_t^{(\alpha)}(\omega) = \inf\left\{y:\int_0^t  \mathbb{1}\left\{Y_s(\omega)\leq y\right\} \rd s > \alpha t\right\}.
\end{equation}
For fixed $\omega$, $M_t^{(\alpha)}(\omega)$ is the $\alpha$-quantile of the function $s\mapsto Y_s(\omega)$, for $s< t$, which is considered as a random variable on the space $([0,t]; \rd s/t)$ equipped with the Borel $\sigma$--field.
\end{Definition}

This definition, whilst it is a quantile diffusion, is different to the one we propose in our construction in the sense that we are not looking directly at the quantiles of a diffusion process.  Rather we produce a quantile diffusion that implicitly induces some underlying diffusion process. This implicit underlying process will not be of primary interest in our formulation, as instead we wish to focus on characterising and constructing diffusions on the quantile space with non--trivial skew--kurtosis and tail characteristics that can be parameterised and interpreted. As such we develop parameterised mappings of an underlying diffusion process in such a way that the resulting quantile diffusions are not the quantiles of the underlying driving diffusion, as produced in Definition \ref{quantileprocessdefn2}, but rather will imply such a process without ever requiring its explicit specification.

A third characterisation of quantile processes is also widely adopted in the statistics and econometrics time series literature, see \cite{koenkerbassett}, where a regression framework is developed for (linear or nonlinear) conditional quantile functions at any, or over all, quantile levels $u\in[0,1]$.  As an example, one may consider the following definition, given in the tutorial review of various models in \cite{quantiletutorial}, which allows for the autoregressive parameters of the model to vary with the quantile level---see also \cite{koenkerqr40,quantilear}.  

\begin{Definition}\label{parametricquantileprocessdefn}
Consider a univariate time series $\{Y_1,\ldots,Y_t,\ldots\}$ for $t\in\mathbb{N}$ and let $\mathscr{F}_t=\sigma(Y_1,\ldots,Y_t)$ denote the natural sigma-algebra of the observed time series.  Let $\bm{\theta}\in\mathbb{R}^d$ be a static vector of model parameters, $u\in[0,1]$ a quantile level, and $\alpha_i(u):[0,1]\rightarrow\mathbb{R}$ be unknown quantile functions, given by Definition \ref{quantilefunctiondef}, for $i=1,\ldots,p$.  The conditional quantile autoregressive QAR(p) model for the conditional quantile function of the random variable $Y_t$, conditioned on the observations of the time series until time $t-1$, is characterised by
\begin{equation}\label{conditionalQARmodel}
Q_{Y_t}(u\vert\mathscr{F}_{t-1};\bm{\theta}) = \sum_{i=1}^p \alpha_i(u) Y_{t-i} + Q_{\epsilon}(u;\gamma),
\end{equation}
where $Q_{\epsilon}(u;\gamma)$ denotes the quantile error function, representing the white noise sequence $(\epsilon_t)_{t=1,2,\ldots}$ with $\gamma\in\mathbb{R}^{d'}$ a vector of static parameters.
\end{Definition}
By construction, Eq. (\ref{conditionalQARmodel}) is a discrete--time, function--valued quantile process for the conditional quantile function of $Y_t$.  This is just one approach to the construction of conditional quantile processes; an alternative non--parametric approach is given in \cite{conditionalqp}.  Additionally, a model of this form also admits a quantile time series model, i.e., the underlying time series model is a sequence of quantile levels. However the direct link between two such models may not always be easily obtained in closed form.  The definition of this model is given as follows for a random--level quantile time series.
\begin{Definition}
Consider the setup given in Definition \ref{parametricquantileprocessdefn}, where $\{y_0,y_1,\ldots,y_{t-1} \}$ are observations of the time series until time $t-1$, and $U_t\sim U[0,1]$ for all $t\in\mathbb{N}$ {i.i.d.}.  The functional time series model with random coefficients in an AR structure is given by 
\begin{equation}\label{functionaltimeseriesmodel}
    y_t=\sum_{i=1}^p\alpha_i\left(U_t \right)y_{t-i}+\alpha_0\left(U_tx \right) = \sum_{i=1}^p\alpha_i\left(U_t \right)y_{t-i}+\epsilon_t
\end{equation}
where 
\begin{equation}
    \sum_{i=1}^p\alpha_i\left(U_t \right)y_{t-i}+\alpha_0(U_t)
\end{equation}
is a monotone increasing function of $U_t$ and $(\epsilon_t)_{t=1,2,\ldots}$ is a white noise sequence, independent of $U_t$.  
\end{Definition}
A framework for the construction of such models with specific properties is detailed in \cite{quantiletutorial}.  Moreover, one may consider the random coefficients in Definition \ref{functionaltimeseriesmodel} to be co--monotonic random functional coefficients, so to define a scalar on a function regression version of Eq. (\ref{conditionalQARmodel}). That is, the AR coefficients are now expressed as monotone functions of a scalar random variable by
\begin{equation}
    Y_t=\sum_{i=1}^p\alpha_i(U_t)Y_{t-i}+Q_{\epsilon}\left(U_t;\gamma\right)
\end{equation}
for i.i.d. $U_t\sim U[0,1]$.  See Example 1 in \cite{quantiletutorial} for more details and a discussion on the advantages of constructing such a model.

The above characterisations of quantile processes is not exhaustive, but represents quantile processes usually encountered in time series statistics and econometrics literature. 

This paper adds a novel class of quantile processes in the diffusion context. We draw attention to the emerging stochastic dynamics of the built continuous--time quantile processes, which are associated to some diffusion process, that is, implied but not necessarily explicitly given. The parameters have direct interpretation on the properties of the implicit process that such a quantile diffusion induces. Perhaps somewhat surprisingly, this will be constructed using an auxiliary driving diffusion process to provide the stochasticity required, but it is important to keep in mind that this driving auxiliary diffusion process is not the induced process upon which the quantile diffusion characterises the law. Herein lies the novelty in our perspective and will ultimately provide the utility of our framework as distinct from the aforementioned approaches. 
\begin{rem}
The definitions in this paper produce a wide range of families of quantile processes, which are both parametric and interpretable, and which give the flexibility to consider non--static distributions. While in the present work the constructed quantile processes are diffusions, we emphasize that the devised method by which these are built, see Sections 2.2 and 2.3, may be used to construct continuous--time quantile processes with jumps. In principle, all that is needed is for the underlying driving process to have jumps.
\end{rem}

\subsection{Quantile processes characterized as solutions to a nested fixed--point problem}\label{fixedpointcharacterisationsection}
We begin by presenting a general system of nested fixed--point equations that will characterise the class of solutions that we seek to develop into a formulation of quantile diffusions.  A nested fixed--point equation is comprised of both an inner fixed--point solution, in composition with an outer fixed--point problem. The solution to the inner fixed--point solution will characterise the behaviour of the solution to the outer fixed--point solution through a mapping that we wish to characterise in a continuous--time process setting.  The key result of interest here is the kind of mapping that preserves the solution type.    
\begin{Definition}\label{fixedpointeqndefinition}
Let $(\Omega,\mathscr{F},\mathbb{P})$ be a probability space and $X:\Omega\rightarrow\mathbb{R}$ be a random variable on the probability space with distribution function $F_{X}(x)=\mathbb{P}(X\leq x):\mathbb{R}\rightarrow [0,1]$. Assume for some $\eta\in\mathbb{R}$ that $F^{\prime\prime}_{X}(x)$ exists and is bounded in the neighbourhood of $\eta$, and that $F^\prime_{X}(\eta)=f_{ X}(\eta)>0$.  It follows that when $\eta$ is the fixed--point solution to
\begin{equation}\label{fixedpointeqn}
F_{X}(\eta)=u
\end{equation}
for fixed $u\in(0,1)$, it represents the unique $u$--quantile of $F_{X}$.
\end{Definition}  

Given a sequence of i.i.d. samples from the distribution $F_{X}$, a representation of, and asymptotic convergence result for, the relationship between the population quantile $\eta$ (i.e., the solution to Eq. (\ref{fixedpointeqn})) and the corresponding sample $u$--quantile is presented in \cite{bahadur} (see \cite{kiefer1967bahadur} for an extension) and given as follows.

\begin{Proposition}\label{fixedpointestimatordefn}
Let $\omega=(X_1,X_2\ldots )$ be a sequence of i.i.d. random variables such that for all $i$, $X_i\sim F_X$.  For any sub--sequence $\omega=(X_1,X_2\ldots X_n)$, for $n=1,2,\ldots$, denote by $Y_n=Y_n(\omega)$ the sample $u$--quantile, and by $Z_n=Z_n(\omega)$ the number of observations $X_i$ in the sample such that $X_i>\eta$. Under the assumptions on $F_X$ given in Definition \ref{fixedpointeqndefinition}, it follows that $$Y_n(\omega)=\eta + \left[\tfrac{1}{n}\left(Z_n(\omega)-n(1-p)\right) f_X(\eta)\right] + R_n(\omega)$$
where $R_n(\omega)=\mathscr{O}(n^{-3/4}\log n)$ as $n\rightarrow\infty$ with probability one.
\end{Proposition}
\begin{proof}
We refer to \cite{bahadur} for a proof.
\end{proof}

We introduce the following ``composite problem'', whereby the equation of interest is composed of an inner and an outer version of Eq. (\ref{fixedpointeqn}) for two distinct distribution functions.

\begin{Definition}\label{compositefixedpointdefn}
Let $(\Omega,\mathscr{F},\mathbb{P})$ be a probability space and $Y,Z:\Omega\rightarrow\mathbb{R}$ be two random variables on the probability space with distribution functions $F_Y$ and $F_Z$, respectively.  For some mapping $m:\mathbb{R}\rightarrow\mathbb{R}$, assume $F_Y,F_Z$ satisfy the differentiability conditions given in Definition \ref{fixedpointeqndefinition} in the neighbourhood of two fixed points $\eta_Y$ and $\eta_Z=m(\eta_Y)$ from each of the distributions, respectively.  The composite problem of the above type corresponds to obtaining the fixed--point solution $\eta_Y$ to the pair of equations
\begin{align}\label{outerfixedpoint}
F_Z\left(m\left(\eta_Y\right) \right) &= u_Z, \\ \label{innerfixedpoint}
F_Y\left(\eta_Y\right) &= u_Y
\end{align}
for fixed $u_Z,u_Y\in(0,1)$.
\end{Definition}

Here, each fixed--point refers to a quantile from the given distribution at some fixed level $u\in(0,1)$.  Considering a problem of this type allows one to modify the inner solution in such a way that the result will have specific, desired distributional properties yet still remain to be a fixed--point solution to a problem of the same type. The quantile levels $u_Y$ and $u_Z$ will be equal in value.  The type of mapping $m$ that preserves the solution type, however, is non--trivial, and we seek to understand the form of the mapping that satisfies the problem in Definition \ref{compositefixedpointdefn}, so to provide a parametric class of such maps that can be used to construct the novel class of quantile processes presented in this paper---see Definition \ref{processdrivendefn}.  Such a map will ensure that $\eta_Z:=m(\eta_Y)$ is a solution to the standalone fixed point problem $F_Z(\eta_Z)=u_Z$.

\begin{rem}\label{qpuniquerem}
The solution to Eq. (\ref{outerfixedpoint}) is not uniquely obtained by any one fixed--point solution to Eq. (\ref{innerfixedpoint}) and mapping $m$.
\end{rem}

The class of functions to which such a map must belong is given by Definition \ref{mappingclassdefn}, and our construction--based approach in Section \ref{quantileprocesses} gives a parametric example of such a map for a wide class of flexible families of distributions.  As this paper focuses on the construction of quantile processes in continuous--time, we extend this problem to the the dynamic setting; that is, we let $(\Omega,\mathscr{F}, (\mathscr{F}_t)_{0\leq t<\infty},\mathbb{P})$ be a filtered probability space and consider two processes $(Y_t)_{0\leq t<\infty}$ and $(Z_t)_{0\leq t<\infty}$ on the space equipped with marginal distribution functions $F_Y(t,y)$ and $F_Z(t,z)$, respectively, at each $t\in(0,\infty)$.  For complete generality, we may also consider a time--dependent mapping $m(t,x)$. The problem in Definition \ref{compositefixedpointdefn},  will now correspond to solving the pair of equations 
\begin{align}\label{outerfixedpoint2}
F_Z\left(t,m\left(t,\eta_t^Y\right) \right) &= u_{t,Z} \\ \label{innerfixedpoint2}
F_Y\left(t,\eta_t^Y\right) &= u_{t,Y}
\end{align}
for each $t\in(0,\infty)$ and fixed $u_{t,Y}, u_{t,Z}\in(0,1)$.

One considers the non--stationary extension of Definition \ref{fixedpointeqndefinition} as follows.  For some $\epsilon>0$ and all $t\in(0+\epsilon,\infty)$, define the time interval $\delta_{t,\epsilon}:=(t-\epsilon,t+\epsilon)$ whereby $\epsilon$ is chosen such that $(Y_t)$ and $(Z_t)$ are marginally locally stationary on $\delta_{t,\epsilon}$, i.e., for all $s\in\delta_{t,\epsilon}$, $F_Y(s,y)=:F_{Y,t}(y):\mathbb{R}\rightarrow\mathbb{R}$ and $F_Z(s,z)=:F_{Z,t}(z):\mathbb{R}\rightarrow\mathbb{R}$. For all $t\in(0+\epsilon,\infty)$,  let $\omega_t=(Y_{t,1},Y_{t,2},\ldots, Y_{t,n})$ be a sequence of i.i.d. random variables such that $Y_{t,i}\sim F_{Y,t}$ for all $i=1,\ldots,n$ and $Y_t^n=Y_t^n(\omega_t)$ is the sequence $u_{t,Y}$--quantile, analogously to Proposition \ref{fixedpointestimatordefn}.  Now, if the mapping $m$ belongs to the space of functions in Definition \ref{mappingclassdefn}, for all $t\in(0+\epsilon,\infty)$, it holds that $Z_t^n=Z_t^n(\omega_t)=m\left(t,Y_t^n(\omega_t)\right)$ will be the sequence $u_{t,Z}$--quantile solution to Eq. (\ref{outerfixedpoint2}), satisfying the convergence result in \cite{bahadur}.  For all $t\in(0+\epsilon,\infty)$ and fixed $u_{t,Z}\in(0,1)$, the fixed--point solution to the time--inhomogeneous equation $F_Z(t,\eta_t^Z)=u_{t,Z}$, given by $\eta_t^Z:=m\left(t,\eta_t^Y\right)$ corresponds to the sequence $u_{t,Z}$--quantile of the distribution $F_{Z,t}(z)$.  As such, if the process $(Z_t)$ is constructed marginally by the transform $Z_t=m(t,Y_t)$ for all $t\in(0+\epsilon,\infty)$, we refer to $(Z_t)$ as a quantile process.

Whilst to motivate the results in this paper we view the characterisation of the dynamic composite problem in the above way, we note that extensions of the result by Bahadur \cite{bahadur} in the case of non--i.i.d. random variables and processes are also given in \cite{dutta1971bahadur,hesse1990bahadur,sen1968asymptotic,sen1972bahadur,wu2005bahadur}.

In summary, the above characterisation allows one to obtain sample quantiles of some dynamic distribution $F_Z(t,z)$ with specific properties, from the sample quantiles of any other dynamic, continuous distribution satisfying the necessary differentiability and stationarity conditions, using a mapping of the class given in Definition \ref{mappingclassdefn}. 

\subsection{Construction I: Random--level quantile diffusions}\label{randomlevelconstruction}
The first type of quantile diffusion is constructed by a composite map applied marginally to a driving process, which in principle could be multi-dimensional, but here will be univariate throughout.  This auxiliary driving process produces the stochasticity of the output quantile diffusion, which characterises a law with relative statistical properties determined by the functions in the map.  Such a construction belongs to the class of problems presented in Section \ref{fixedpointcharacterisationsection}, and provides a parameterisation of the type of mapping $m$ given in Definition \ref{mappingclassdefn}. Whilst the parameterisation is simple, it produces a wide class of models with appropriately chosen attributes from a modelling perspective. We motivate this construction by seeking models with some desired statistical properties, both marginally and serially, and we wish to observe the behaviour of the quantiles, or the tail--behaviour, of this model.  The output quantile process will characterise this information without the need to explicitly state the distributional form of the model itself. 
 
\begin{Definition}\label{processdrivendefn}
Let $Q_{\zeta}(u;\bm{\xi})$ be the quantile function of some real--valued random variable $\zeta$, as specified in Definition \ref{quantilefunctiondef} for quantile level $u\in[0,1]$ and $\bm{\xi}\in\mathbb{R}^d$ a $d$--dimensional vector of constant parameters.  Consider a process $(Y_t)_{0\leq t<\infty}$ specified by Definition \ref{diffusiondef}.  Two cases are considered next:
\begin{enumerate}[(i)] 
\item At any time $t\in[t_0,\infty)$, $t_0>0$, the process $(Y_t)$ is governed by a continuous, time--dependent distribution function $F_Y(t,y; \bm{\theta}):\mathbb{R}^+\times\mathbb{R}\rightarrow[0,1]$ where $\bm{\theta}\in\mathbb{R}^{k}$ is a $k$--dimensional vector of constant parameters, such that the marginals of the process given by $U_t:=F_Y(t,Y_t{;\bm{\theta}})$ are {\it uniformly distributed} on $[0,1]$. We say that $F_Y(t,y;\bm{\theta})_{t_0\leq t<\infty}$ is the ``{\it true law}'' of $(Y_t)$. 

\item At any time $t\in[t_0,\infty)$, let $F(t,y;\bm{\vartheta}):\mathbb{R}^+\times\mathbb{R}\rightarrow[0,1]$ be a continuous, time--dependent distribution function, with $\bm{\vartheta}\in\mathbb{R}^{k^{\prime}}$ a $k^{\prime}$--dimensional vector of constant parameters, such that the process given by $\widetilde{U}_t:=F(t,Y_t;\bm{\vartheta})=F(t,Q_Y({t,} U_t ;\bm{\theta});\bm{\vartheta})$ is {\it non-uniformly distributed} on $[0,1]$. We say that $F(t,y;\bm{\vartheta})_{t_0\leq t<\infty}$ is the ``{\it false law}'' of $(Y_t)$.
 \end{enumerate}
At each time $t\in[t_0,\infty)$, the random--level quantile diffusion is defined by
\begin{equation}\label{processdrivenQP}
Z_t = Q_{\zeta}\left(U_t;\bm{\xi}\right),
\end{equation}
that is $Z_t=Q_{\zeta}(F_Y(t,Y(t,\omega);\bm{\theta});\bm{\xi}):\mathbb{R}^+\times(\mathbb{R}^+\times\Omega)\rightarrow[-\infty,\infty]$ in case (i). In case (ii), the quantile diffusion is defined by
\begin{equation}\label{processdrivenQP2}
Z_t = Q_{\zeta}\left(\widetilde{U}_t;\bm{\xi}\right) = Q_{\zeta}\left(F\left(t,Q_Y\left(t,U_t;\bm{\theta}\right);\bm{\vartheta}\right);\bm{\xi}\right),
\end{equation}
that is, $Z_t=Q_{\zeta}(F(t,Y(t,\omega);\bm{\vartheta});\bm{\xi}):\mathbb{R}^+\times(\mathbb{R}^+\times\Omega)\rightarrow[-\infty,\infty]$. In each case, $Z_t:[0,1]\rightarrow[-\infty,\infty]$ is an output to the input $U_t$ for all $t\in[t_0,\infty)$.  Here, $\zeta$ characterises the family of quantile diffusions to which $(Z_t)$ belongs.
\end{Definition}

The process $(Z_t)$ is well--defined for any choice of initial time $t_0=0+\epsilon$, $\epsilon>0$, as we can ensure the marginal distribution $F_Y(t,y;\bm{\theta})$ is continuous for any $t>0$. At $t=0$, the probability mass is concentrated on $y_0\in\mathbb{R}$.  We note that $F_Y(t,y;\bm{\theta})$ depends on $y_0\in\mathbb{R}$ (as a parameter) so for each $y_0$ we get a different quantile process $(Z_t)$, each of which will satisfy the Markov property. In principle, one could also set $Z_0=z_0\in \mathbb{R}$ and extend the time interval, on which $(Z_t)$ is defined, to $t\in[0,\infty)$.

In the two cases in Definition \ref{processdrivendefn}, at each $t\in[t_0,\infty)$, we view the value of the process $(U_t)$ as the quantile level that each marginal of the output quantile diffusion corresponds to. Since $(U_t)$ is uniformly distributed on $[0,1]$ at each $t\in[t_0,\infty)$, the output quantile process yields well--defined quantiles for all quantile levels in $[0,1]$.  Unlike the quantile process definitions given in Section \ref{quantileprocesses}, the quantile level that we observe is itself a random process that evolves over time on the interval $[0,1]$, driving the change in quantile value given by the output process $(Z_t)$.

We emphasise that in case (ii), the map from the driving process $(Y_t)$ to the uniformly distributed process $(U_t)$ has two stages: first, we map to $\widetilde{U}_t=F(t,Y_t;\bm{\vartheta})$, and then we map to $U_t=Q_Y(t,\widetilde{U}_t;\bm{\theta})$ for each $t\in[t_0,\infty)$.  Whilst in general $F_Y$ and $Q_Y$ may be unknown in this case, it is feasible to estimate them.

The process $(U_t)$ is referred to as a ``uniformized diffusion process'' in \cite{bibbona}, and its dynamics are derived.  It is also stated that the same $(U_t)$ may be constructed from different driving processes and their marginal distributions, however this is equivalent to the driving processes having the same serial--dependence as characterised by a unique copula, see \cite{nelsencopula}, on the Fokker--Planck transition distribution.  We may use the well--known result, given in \cite{nelsencopula} by an application of Sklar's theorem \cite{sklar}, that the copula of continuous random variables is invariant under a monotonic transformation of those random variables to highlight that the copula of the transition distribution of the driving process characterises that of the transition distribution of the output quantile process.  Here, our focus lies in the distortion of the marginal distribution of the process $(Y_t)$ to produce $(Z_t)$.  Whilst marginally distinct constructions of quantile diffusion may admit the same dynamic distributional representation, the joint distribution functions of the marginals of the process at subsequent time points will allow for distinct statistical properties between the processes.  The driving process is chosen to establish these desired properties in $(Z_t)$.

One can interpret the resulting quantile process in two ways:  (i) For a fixed time $t\in[t_0,\infty)$, and for each realisation $Y(t,\omega)$ of the underlying driving process $(Y_t)$, the random--level quantile diffusion, defined by either Eq. (\ref{processdrivenQP}) or (\ref{processdrivenQP2}), is scalar--valued and corresponds to a single, fixed quantile level. (ii) The quantile diffusion may be viewed from a path--based perspective, where we observe scalar--valued sequences of quantiles corresponding to some sequence of quantile levels as time evolves.

\subsection{Construction II: Function--valued quantile diffusions}\label{functionvaluedconstruction}
We now define an alternative way to construct stochastic quantile processes, by introducing  the {\it function--valued} quantile diffusion.  This construction leads to diffusions with dynamic statistical properties, which are captured in the parameters of the chosen quantile function. 

This quantile diffusion is a function--valued process; for a quantile function $Q_{\zeta}$ used in the construction, the function space is that characterised by the corresponding distribution function $F_{\zeta}=Q_{\zeta}^{-}$.

\begin{Definition}\label{parameterdrivendefn}
Consider $(\Omega,\mathscr{F},(\mathscr{F}_t)_{0\leq t<\infty}, \mathbb{P})$ and the $(\F_t)$-adapted, $d$-dimensional process $(\bm{\xi}_t)_{0\leq t<\infty}$, satisfying the multivariate version of the SDE (\ref{ysde}) with $\bm{\xi}_0\in\mathbb{R}^d$, where the Brownian motion $(W_t)_{0\leq t<\infty}$ is $n$-dimensional (not necessarily $n=d$) with $\rd W_t^{(i)}\rd W_t^{(j)}=\rho_{ij}\rd t$, $\rho_{ij}\in[-1,1]$ for $i\neq j$. Furthermore, let $(\bm{\xi}_t)$ be a vector of stochastic parameters and define $Q_{\zeta}(u;\bm{\xi}_t)$ such that for each $\omega\in\Omega$, $t\in[0,\infty)$, $\bm{\xi}_t=\bm{\xi}(t,\omega)\in\mathbb{R}^d$ and $Q_{\zeta}(u;\bm{\xi}(t,\omega))$ is a quantile function given by Definition \ref{quantilefunctiondef}. Then, for $t\in[0,\infty)$, the function--valued quantile diffusion is defined by
\begin{equation}\label{parameterdrivenQP}
Z_t(u) = Q_{\zeta}\left(u;\bm{\xi}_t\right),
\end{equation}
where $u\in[0,1]$ is the quantile level.
We have $Z_t(u)=Z(t,\omega,u):\mathbb{R}^+_0\times\Omega\times [0,1]\rightarrow [-\infty,\infty]$. For each $\omega\in\Omega$ and $t\in[0,\infty)$, $Z(t,\omega,u)=Q_{\zeta}(u;\bm{\xi}(t,\omega))\in \mathscr{C}^2\left([0,1] \right)$ in the case that $F_{\zeta}(z)\in\mathscr{C}^2\left(\mathbb{R}\right)$ where $\mathscr{C}^2$ is the space of continuously twice differentiable functions. 
\end{Definition}

One can relax the differentiability requirements in the aforementioned Definition \ref{parameterdrivendefn} with regards to the space to which $F_{\zeta}(z)$ belongs. 

For each fixed $t$ and $\omega$, $Z(t,\omega,u)$ are functions of the quantile level $u$ and hence, by construction, are elements from the space of quantile functions $\mathbb{Q}_{\zeta}$. The function--valued quantile process takes values in a function--valued space because each realisation $\bm{\xi}(t,\omega)$ at some time $t\in[0,\infty)$, yields the value of the quantile process over all quantile levels $u$. 
To obtain a quantile process similar to that in Definition \ref{quantileprocessdefn}, one could construct a quantile diffusion of this type and fix the quantile level at some $\bar{u}\in[0,1]$.  We treat this special case, next.

\subsection{Link between the two constructions of quantile diffusions}
In some unusual cases, one may observe a direct connection between the random--level and function--valued constructions.  We consider the following case in which we can construct a version of random--level quantile diffusions given in Definition \ref{processdrivendefn}, however now where the underlying driver is a stochastic vector of parameters ($\bm{\xi}_t)$, as in Definition \ref{parameterdrivendefn}, and we have control over the fixed quantile level corresponding to the quantiles modelled by the output process. 
\begin{Definition}\label{parameterprocessdrivendefn}
Let $(\bm{\xi}_t)$ be the stochastic vector of parameters given in Definition \ref{parameterdrivendefn}, $\bar{u}\in[0,1]$ be a fixed quantile level, and $Q_{\zeta}(u;\bm{\xi}_t)$ be defined as per Definition \ref{parameterdrivendefn}.  Consider the special case of the function--valued construction in Definition \ref{parameterdrivendefn}, whereby we fix $u=\bar{u}\in[0,1]$, and so the quantile process given by Eq. (\ref{parameterdrivenQP}) becomes
\begin{equation}\label{ztubar}
Z_t^{\bar{u}}:=Z_t(\bar{u}) =Q_{\zeta}\left(\bar{u};\bm{\xi}_t\right),
\end{equation}
where $Z_t^{\bar{u}}=Q_{\zeta}(\bar{u};\bm{\xi}(t,\omega)):\mathbb{R}^+\times\Omega\rightarrow[-\infty,\infty]$.  This is distinct from the usual case whereby the function--valued construction models the dynamics of the entire quantile curve, and is equivalent to taking some fixed point on the quantile curve in the function--valued construction.  

Now, let the functions $Q_{\tilde{\zeta}}(u;\bm{\tilde{\xi}})$ and $F(t,y;\bm{\vartheta})$ be the quantile function and distribution function, respectively, given in Definition \ref{processdrivendefn} (ii). The process analogous to the random--level quantile diffusion, however now at fixed level $\bar{u}$ and with the stochastic driver $(\bm{\xi}_t)$, is defined by
\begin{equation}\label{parameterprocessdrivenQP2}
Z_t = Q_{\tilde{\zeta}}\left(F\left(t,Z_t^{\bar{u}};\bm{\vartheta}\right);\bm{\tilde{\xi}}\right)
\end{equation}
where $Z_t^{\bar{u}}$ is given by Eq. (\ref{ztubar}). Since Eq. (\ref{parameterprocessdrivenQP2}) is implicitly dependent on $\bar{u}$, this process takes well--defined quantile values at the quantile level $\bar{u}$.         
\end{Definition}

If one chooses the functions $Q_Y$ and $Q$ in Definition \ref{processdrivendefn} (ii) such that $Q_Y(t,U_t)\overset{d}{=}Z_t^{\bar{u}}$, where $\overset{d}{=}$ denotes equal in distribution, for each $t\in[t_0,\infty)$, $t_0>0$, one can ensure this quantile diffusion matches that obtained by the usual random--level construction given in Eq. (\ref{processdrivenQP2}).

\subsubsection{Example I}
Consider a uniformly distributed random variable $X\sim U[a,b]$ where $-\infty<a<b<\infty$.  Take $a=0$ and let the parameter $(b_t)_{0\leq t<\infty}$ be stochastic, satisfying the SDE \ref{ysde} where $\mu(t,x):\mathbb{R}^+\times\mathbb{R}^+\rightarrow\mathbb{R}$ and $\sigma(t,x):\mathbb{R}^+\times\mathbb{R}^+\rightarrow\mathbb{R}^+$ satisfy the necessary conditions to ensure $b_t>0$ for all $t\in[0,\infty)$.  Using Definition \ref{parameterdrivendefn}, we construct a uniformly--distributed function--valued quantile diffusion by $Z_t(u) = Q_X(u;b_t) = ub_t$ for $u\in[0,1]$.  Fixing $u=\bar{u}\in[0,1]$, the process $Z_t^{\bar{u}}=\bar{u}b_t$ satisfies the SDE $\rd Z_t^{\bar{u}} =\bar{u}\mu(t,Z_t^{\bar{u}}/\bar{u})\rd t+ \bar{u}\sigma(t,Z_t^{\bar{u}}/\bar{u})\rd W_t$.

Taking this to be the SDE satisfied by the driving process $(Y_t)$ in Definition \ref{processdrivendefn}, one obtains a special case of a random--level quantile diffusion driven by the parameter process $(b_t)$, producing output quantiles at level $\bar{u}\in[0,1]$.  

\subsubsection{Example II}
Consider some random variable $X_2\sim F_{X_2}$ that belongs to the location-scale family with location parameter $A\in\mathbb{R}$ and scale parameter $B\in\mathbb{R}^+$, that is $X_2\overset{d}{=}A+BX_1$ for any random variable $X_1\sim F_{X_1}$.  Take $B=1$ and let the location parameter $(A_t)_{0\leq t<\infty}$ be stochastic, satisfying the SDE \ref{ysde} with associated law $F_A(t,a)_{0<t<\infty}$.  Using Definition \ref{parameterdrivendefn}, we construct a location--scale, function--valued quantile diffusion by $Z_t(u) = Q_{X_1}(u) +A_t$ for all $t\in[0,\infty)$ and $u\in[0,1]$.  Fix $u=\bar{u}\in[0,1]$ and define a distribution function $F_Y(t,y)=F_A(t,y-Q_{X_2}(\bar{u}))$ for all $t\in(0,\infty)$. For some choice of the functions $Q_{\zeta}$ and $F$ in Definition \ref{processdrivendefn}, one can produce equivalent quantile diffusions in the two following ways:
\begin{enumerate}
\item Using the function--valued construction, taking Eq. (\ref{processdrivenQP2}) with $Q_Y$ the quantile function corresponding to the distribution function $F_Y(t,y)=F_A(t,y-Q_{X_2}(\bar{u}))$.
\item By $Z_t = Q_{\zeta}(F(t,Z_t^{\bar{u}}))$ where $Z_t^{\bar{u}}=Q_{X_1}(\bar{u}) +A_t$ for each $t\in(0,\infty)$. This is a special case of a random--level quantile diffusion, where the driving process is the location parameter process $(A_t)$ and $(Z_t)$ models quantiles at the chosen level $\bar{u}$.
\end{enumerate} 

We emphasise that the above examples illustrate the ``special case'' of quantile diffusions that is not a new construction but instead an overlap between the random--level and function--valued constructions.  In general, however, these two constructions will not be related.  The main difference between them is that, at any fixed time, the random--level quantile diffusions lie in Euclidean space and the function--valued ones in a function space.  We emphasise that, when one observes the output processes of each of these constructions from a path--based perspective, the random--level construction presents as a stochastic sequence of well--defined quantiles (each of which corresponds a random quantile level that is determined by either $F_Y$ or $F$, and the driving process).  The function--valued construction presents as a sequence of quantile curves from the same family but each with different parameters (which are determined by the values of the driving parameter diffusion).  In what follows, we focus on the random--level construction.

\subsection{Quantile transforms: rank transmutation map}\label{familiesofquantiletransforms}
In this section we discuss types of transformations that have been developed and allow one to map from the distributional setting to the quantile setting with flexibility in the properties of the quantiles produced.  Rank transmutation maps (RTMs), as introduced in \cite{rtm}, are a composition of a base quantile function and some target distribution function and enable one to produce more flexible classes of quantile functions, relative to the base.  Such maps are also referred to as $P$-transformations in \cite{gilchrist}.  

\begin{Definition}\label{rtmsdefn}
Consider two distribution functions $F_1$ and $F_2$ with a common sample space.  A pair of general RTMs may be given by
\begin{align}
&G_{R_{12}}(u) = F_2(F_1^{-}(u)),&
&G_{R_{21}}(u) = F_1(F_2^{-}(u)),&
\end{align}
where each of these functions maps $[0,1]$ onto itself. Under suitable assumptions, $G_{R_{12}}(u)$ and $G_{R_{21}}(u)$ are mutual inverses that satisfy $G_{R_{ij}}(0)=0$ and $G_{R_{ij}}(1)=1$, for $i,j=1,2$ and $i\neq j$.  
\end{Definition}
The assumption that the RTMs be continuously differentiable is made in order to ensure that the densities of the mapped random variables are continuous, and one may also assume that they be monotone so that the distribution and quantile functions involved are well--defined.  Families of RTMs, including but not limited to a quadratic class, skew--uniform, skew--exponential and skew--kurtotic classes, are outlined in \cite{rtm}.  Each of these classes allow to devise different properties of the distorted quantiles, and are discussed in \cite{quantiletutorial} in the context of quantile time series models. Another family of quantile distortion maps is the so--called Tukey elongation transforms, as detailed in \cite{peterssisson,lmomentsgh}, where the idea is to construct skewed or heavy--tailed distributions by transforming some base random variable, which is often taken to be Gaussian.  The amount of skewness or kurtosis introduced is relative to the base random variable.  In Section \ref{familiesofmodelssection}, we focus on the $g$--transform, $h$--transform and $g$--$h$ subfamilies; the $g$--$k$ subfamily has also been largely studied in the literature.  
\begin{Definition}\label{elongationmapdef}
An elongation transform $T:\mathbb{R}\rightarrow\mathbb{R}$ applied to a generic random variable $W$ is a map that satisfies $T(w)=T(-w)$, $T(w)=w+\mathscr{O}(w^2)$ for $w$ around the mode, and is such that $T^{\prime}(w)>0$ and $T^{\prime\prime}(w)>0$ for $w>0$.
\end{Definition}  
\begin{Definition}\label{tukeytransform-general}
Consider some base random variable $W\sim F_W$ and the transformation $X := WT(W)^\theta$ for a parameter $\theta\in\mathbb{R}$ and $T:\mathbb{R}\rightarrow\mathbb{R}^+$ a Tukey elongation map satisfying Definition \ref{elongationmapdef}.  The random variable $X$ is Tukey--distributed and has quantile function $Q_X(u) = A + BQ_W(u)T(Q_W(u))^\theta$.  For $u\in[0,1]$, $Q_W(u)$ is the quantile function of the base random variable; $A\in\mathbb{R}$ and $B\in\mathbb{R}^+$ are the location and scale parameters, respectively.  
\end{Definition}
In order to generate a more pronounced relative kurtosis when compared to the base random variable, one can use the $h$-- and $k$--transform classes where $T(W)$ is given by $T_h(w)=\expp(w^2)$ or $T_k(w)=1+w^2$, respectively.  The $g$--transform class, where $T(W)$ is given by $T_g(w)=(\expp(w)-1)/w$ can be used to introduce relative skewness.  Further details on the Tukey elongation subfamilies and their application to quantile error functions with desirable properties are given in \cite{quantiletutorial}.  These mappings may be seen as a parametric family of RTMs if in Definition \ref{rtmsdefn} one has $F_1 = F_W$ such that $W = F_1^{-}(U)$ for $U\sim U[0,1]$ and $T(W) = G_{R_{12}}(U) = F_2(F_1^{-}(U))$. At each instance in time, these quantile transforms are analogous to the random--level quantile diffusion map in Definition \ref{processdrivendefn}.  The random--level case produces a dynamical evolution of these distortion maps as we move from the setting where we transform some base random variable to transforming a base stochastic process.

\section{Quantile diffusion SDEs} \label{QuantDiffSDE}
In this section we derive the dynamical equations of the quantile diffusions. We begin by deriving SDEs for the dynamics of the random--level quantile diffusion, considering the two cases in Definition \ref{processdrivendefn}. The following proposition and corollary guarantee that in each case $(Z_t)$ is a diffusion, and provide the corresponding infinitesimal drift and volatility coefficients.  As such, the transition distribution of $(Z_t)$ must solve the Kolmogorov backward equation having these same coefficients. We omit the dependence on the vectors of parameters in the notation for the following distribution, quantile and density functions.

\begin{Proposition}\label{prop1} 
Let $(Z_t)_{t_0\leq t<\infty}$ be a quantile diffusion given by Definition \ref{processdrivendefn} (ii). Assume the following derivatives exist so that $f(t,y):=\partial_y F(t,y)$ and $f_{\zeta}(z):=\partial_z F_{\zeta}(z)$ are density functions. The dynamics of $(Z_t)$ satisfies
\begin{equation}\label{dz1}
\rd Z_t=\alpha(t,Z_t) \rd t+\widetilde{\sigma}(t,Z_t) \rd W_t
\end{equation} 
where
\begin{align}\label{drift}\begin{split}
\alpha(t,Z_t) &= \frac{\partial_t F(t,x)\vert_{x=Q(t,F_{\zeta}(Z_t))}}{f_{\zeta}(Z_t)}+\mu(t,Q(t,F_{\zeta}(Z_t)))\frac{f(t,Q(t,F_{\zeta}(Z_t)))}{f_{\zeta}(Z_t)} 
\\[0pt] &\qquad+\frac{1}{2}\sigma^2(t,Q(t,F_{\zeta}(Z_t)))\frac{f'(t,Q(t,F_{\zeta}(Z_t)))f_{\zeta}(Z_t)^2 - f(t,Q(t,F_{\zeta}(Z_t)))^2f_{\zeta}'(Z_t)}{f_{\zeta}(Z_t)^3}, \end{split} 
\\[0pt] \label{volatility} \widetilde{\sigma}(t,Z_t) &= \sigma(t,Q(t,F_{\zeta}(Z_t)))\frac{f(t,Q(t,F_{\zeta}(Z_t)))}{f_{\zeta}(Z_t)},
\end{align}
for $t\in[t_0,\infty)$ and $Z_{t_0}=z_{t_0}\in\mathbb{R}$. The short-hand notation $f'$ denotes differentiation with respect to the spatial variable.
\end{Proposition}
\begin{proof}
The result follows from a straightforward application of Ito's formula and also by use of Eq. (64) in \cite{quantilemech}.
\end{proof}

\begin{Corollary}\label{corr1}
Let $(Z_t)$ be a quantile diffusion given by Definition \ref{processdrivendefn} (i). The dynamics of $(Z_t)$ satisfy
\begin{equation}\label{dz2}
\rd Z_t = \alpha(t,Z_t)\rd t+\widetilde{\sigma}(t,Z_t)\rd W_t
\end{equation}
where
\begin{align}
\begin{split}
\alpha(t,Z_t) &= \frac{\sigma^2(t,Q_Y(t,F_{\zeta}(Z_t)))f_Y'(t,Q_Y(t,F_{\zeta}(Z_t)))}{f_{\zeta}(Z_t)} + \frac{f_Y(t,Q_Y(t,F_{\zeta}(Z_t)))(\sigma^{2}(t,Q_Y(t,F_{\zeta}(Z_t))))^{\prime}}{2f_{\zeta}(Z_t)} \\[0pt] &\qquad-\frac{1}{2}\sigma^2(t,Q_Y(t,F_{\zeta}(Z_t)))\frac{f_Y(t,Q_Y(t,F_{\zeta}(Z_t)))^2f_{\zeta}^{\prime}(Z_t)}{f_{\zeta}(Z_t)^3},
\end{split} \\[0pt]
\widetilde{\sigma}(t,Z_t) &=\sigma(t,Q_Y(t,F_{\zeta}(Z_t)))\frac{f_Y(t,Q_Y(t,F_{\zeta}(Z_t)))}{f_{\zeta}(Z_t)},
\end{align}
for $t\in[t_0,\infty)$ and $Z_{t_0}=z_{t_0}\in\mathbb{R}$, where $f_Y(t,y)$ is the transition density of the driving process $(Y_t)$ starting with $y_0\in\mathbb{R}$.
\end{Corollary}

\begin{proof}
Similarly to the proof of Proposition (\ref{prop1}), we apply It$\hat{\rm{o}}$'s formula to $Z_t = Q_{\zeta}(F_Y(t,Y_t))$.  Since $F_Y(t,x)$ is the law of the process $(Y_t)$, we can use the Fokker-Plank equation to describe how the density of $(Y_t)$, that is $f_Y(t,y)$, evolves with time. The chain rule yields $\partial_t Q_{\zeta}(F_Y(t,y)) = \partial_t F_Y(t,y)/f_{\zeta}(Q_{\zeta}(F_Y(t,y))$
and by the fundamental theorem of calculus, we obtain
\begin{equation*}
   \partial_t \left( \int_{-\infty}^{\phi(t)}f_Y(t,x)dx\right) = f_Y(t,\phi(t))\partial_t\phi(t)+\int_{-\infty}^{\phi(t)}\partial_t f_Y(t,x)dx.
\end{equation*}
Now, using the Fokker-Planck equation for the marginal density of $(Y_t)$, we have \\ $\partial_t F_Y(t,y) = \int_{-\infty}^y \partial_t f_Y(t,x)dx= -\mu(t,y)f_Y(t,y) +\frac{1}{2}\left(\sigma^2(t,y)f_Y^{\prime}(t,y)+f_Y(t,y)\partial_y\sigma^2(t,y)\right)$, and therefore 
\begin{equation*}
\partial_t \left(Q_{\zeta}\left(F_Y(t,y)\right)\right) = \frac{-\mu(t,y)f_Y(t,y)+\frac{1}{2}\left(\sigma^2(t,y)f_Y^{\prime}(t,y)+f_Y(t,y)\partial_y\sigma^2(t,y) \right)}{f_{\zeta}(t,Q_{\zeta}(t,F_Y(t,y)))}.
\end{equation*}
Noting that $Y_t = Q_Y(t,F_{\zeta}(Z_t))$, the result stated in the corollary follows.
\end{proof} 

We highlight that the choice of the distribution function $F$ in the composite map impacts the drift function $\alpha(t,z)$ of the resulting quantile diffusion SDE.  In the case where $F=F_Y$ is the ``true law'' of the driving process, any explicit dependence of $\alpha(t,z)$ on the drift function $\mu(t,y)$ of the driving process is removed.

Similarly, the SDE satisfied by the function--valued quantile diffusion in Definition \ref{parameterdrivendefn} is obtained by a straightforward application of Ito's formula as follows.  

Let $Z_t(u)_{0\leq t<\infty}$ be a quantile diffusion given by Definition \ref{parameterdrivendefn}.  For each $u\in[0,1]$, the dynamics of $Z_t(u)$ are given by the SDE
\begin{equation}\begin{split}
\rd Z_t(u) &= \sum_{i=1}^d\left[\frac{\partial Q_{\zeta}(u;\bm{\xi}_t)}{\partial\xi_t^{(i)}}a_i\left(t,\xi_t^{(i)}\right) + \frac{1}{2}\sum_{j=1}^d\frac{\partial^2 Q_{\zeta}(u;\bm{\xi}_t)}{\partial\xi_t^{(i)}\xi_t^{(j)}}b_i\left(t,\xi_t^{(i)}\right)b_j\left(t,\xi_t^{(j)}\right)\rho_{ij}\right] \rd t \\[0pt] &\qquad+ \sum_{i=1}^d\frac{\partial Q_{\zeta}(u;\bm{\xi}_t)}{\partial\xi_t^{(i)}}b_i\left(t,\xi_t^{(i)}\right)\rd W_t^{(i)}.
\end{split}\end{equation}

A discussion on the existence of strong and weak solutions to SDEs is given in Appendix \ref{strongsolutionsection}

\section{Tukey quantile diffusions}\label{familiesofmodelssection}
In this section, we develop the Tukey one--parameter $g$--transform and $h$--transform families, and the more general and flexible two--parameter Tukey $g$--$h$ family, of quantile diffusions.  These distributions, which are defined by their quantile functions, were first introduced in \cite{tukey} and allow for flexible modelling of asymmetry and leptokurtosis. 

Let $X$ be a continuous, symmetric random variable that is normalised to have mean zero and variance one.  Throughout this section we consider the standard normal random variable $X\sim\mathscr{N}(0,1)$.

\subsection{Tukey--$g$ transform quantile diffusions}\label{gtransformsection}
\begin{figure}[H]
\setlength\fwidth{.75\textwidth}
\centering
\input{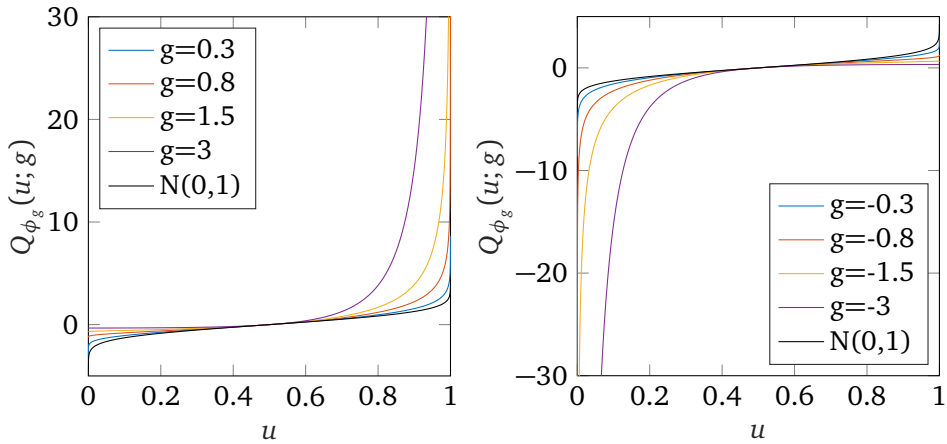}
\caption{$g$--transform quantile functions for $g\in\{0.3,0.8,1.5,3\}$ and $g\in\{-0.3,-0.8,-1.5,-3\}$ relative to the standard normal quantile function.}
\label{gtransformquantile}
\end{figure}

We begin by considering the $g$--transform class of distributions, where a random variable with this distribution is generated by a transformation of the random variable $X$. Denote the Tukey--$g$ transformation map, belonging to the class of transforms presented in Definition \ref{tukeytransform-general}, by $T_{A,B,g}(x)$. The random variable $\phi_g$, the distribution of which will fall under the $g$--transform class, is defined by
\begin{equation}
\phi_g := T_{A,B,g}(X) := A+\frac{B}{g}\left(\expp\left(gX\right)-1 \right),
\end{equation} 
where $A\in\mathbb{R}, B\in\mathbb{R}^+$ are the location and scale parameters, respectively, and the parameter $g\in\mathbb{R}\setminus \{0\}$ allows to flexibly model the skewness. We consider the single--parameter $g$--transform family here, however more flexibility can be introduced by allowing for the skew parameter to have a polynomial representation.  
\begin{Definition}\label{gtransformquantilefunction}
The quantile function of the single parameter $g$--transform family of distributions is given by $Q_{\phi_g}(u;\bm{\xi}) = T_{A,B,g}(x_u)$
for quantile level $u\in[0,1]$, where $x_u=\sqrt{2}\erf^{-}(2u-1)$ is the $u^{th}$ quantile of the standard normal distribution and $\bm{\xi}=(A,B,g)$ is a vector of parameters.
\end{Definition}
There are two cases: $g\in(0,\infty)$ and $g\in(-\infty,0)$, where we have $Q_{{\phi_g}}(u;\bm{\xi}): [0,1]\longrightarrow \left[A-B/g,\infty\right)$ and $Q_{\phi_g}(u;\bm{\xi}): [0,1]\longrightarrow \left(-\infty, A -B/g\right]$, respectively,
excluding the cases where $g=\pm\infty$, i.e., a perfectly skewed model, and $g=0$, i.e., no skewness is introduced.  Figure \ref{gtransformquantile} shows the $g$--transform quantile function, for different values of the parameter $g$, relative to the quantile function $x_u$ of a standard normal random variable.  
\begin{Definition}\label{gtransformqp}
The random--level $g$--transform quantile diffusion is given by Definition \ref{processdrivendefn} where we take $Q_{\zeta}(u;\bm{\xi})=Q_{\phi_g}(u;A,B,g)$. Here, $(Z_t)_{t_0\leq t<\infty}$ has the form 
\begin{equation*}
Z_t=A+\frac{B}{g}(\expp(g\sqrt{2}\erf^{-}(2\,F_{\bullet}(t,Y_t ;\bm{\theta})-1))-1),
\end{equation*}
for all $t\in[t_0,\infty)$ where $F_{\bullet}(t,y ;\bm{\theta})$ may or may not be the law of the process $(Y_t)$.
\end{Definition}

Proposition \ref{prop1} can be used to obtain the dynamics of the $g$--transform quantile diffusion, as shown in \ref{gtransformcoefficients}.  The drift and diffusions coefficients derived next will ensure a unique solution to the SDE \ref{dz1} exists if they satisfy the conditions discussed in Appendix \ref{strongsolutionsection}.

\subsection{SDE coefficients of the $g$--transform quantile diffusion}\label{gtransformcoefficients}
The dynamics of the $g$--transform quantile diffusion are obtained by use of Proposition \ref{prop1}, where we refer to \cite{ghdensity} for the density function $f_{\phi_g}$ of the $g$--transform distribution, given by $f_{\phi_g}(z)=\expp\left(-0.5x_u^2 \right)/ (\sqrt{2\pi}B\expp\left(gx_u \right))$ with $u=F_{\phi_g}(z)$ and where $\bm{\xi}=(A,B,g)$.  
Using the resulting expressions, taking $A=0, B=1$ with no loss of generality for the standardised case, we have the following for the drift function of the $g$--transform quantile diffusion:
\begin{equation}\label{gtransformdrift1}\begin{split}
\alpha(t,Z_t) =&\left\{\mu\left(\cdot,\cdot\right) \widetilde{f}\left(Z_t\right) + \frac{1}{2}\sigma^2\left(\cdot,\cdot\right)  \widetilde{f}'\left(Z_t\right)\right\}\sqrt{2\pi}(gZ_t+1)\expp\left(\frac{\left(\lnn(gZ_t+1\right)^2}{2g^2} \right)   \\[0pt] &+ \sigma^2\left(\cdot,\cdot\right) \widetilde{f}\left(Z_t\right)^2\pi(gZ_t+1)\left\{g + \frac{\lnn(gZ_t+1)}{g} \right\} \expp\left(\frac{\left(\lnn(gZ_t+1)\right)^2}{g^2} \right)
\end{split}\end{equation} 
for a time--homogeneous $F$ that is the ``false law'' of the driving process, and 
\begin{equation}\label{gtransformdrift2}\begin{split}
\alpha(t,Z_t) =&\sigma^2\left(\cdot,\cdot\right) \widetilde{f}^{\prime}\left(Z_t\right)\sqrt{2\pi}\left(gZ_t+1\right)\expp\left(\frac{\left(\lnn(gZ_t+1\right)^2}{2g^2} \right) 
\\[0pt] 
&+ \sigma^2\left(\cdot,\cdot\right) \widetilde{f}\left(Z_t\right)^2 \pi(gZ_t+1)\left\{g + \frac{\lnn(gZ_t+1)}{g} \right\}\expp\left(\frac{\left(\lnn(gZ_t+1)\right)^2}{g^2} \right)
\end{split}\end{equation} 
when $F=F_Y$ is the `true law' of the process $(Y_t)$ with $f_Y(t,y)$ the corresponding transition density.  The volatility function is given by 
\begin{equation}\label{gtransformvol}
\widetilde{\sigma}(t,Z_t) = \sigma\left(\cdot,\cdot\right) \widetilde{f}\left(Z_t \right)\sqrt{2\pi}\left(gZ_t+1\right)\expp\left(\frac{\left(\lnn (gZ_t+1)\right)^2}{2g^2} \right).
\end{equation}    

The arguments of the drift and diffusion parameters are given by $\mu(\cdot,\cdot):= \mu(t, \widetilde{Q}(t,Z_t))$ and $\sigma(\cdot,\cdot):= \sigma(t, \widetilde{Q}(t,Z_t))$, respectively. We define $\widetilde{Q}\left(t,Z_t\right) = Q(t,0.5( 1+\erf( \lnn(gZ_t + 1)/g\sqrt{2})))$, $\widetilde{f}\left(Z_t\right) = f\left(t,\widetilde{Q}\left(t,Z_t\right)\right)$, $\widetilde{f}^{\prime}\left(Z_t\right) = \partial \widetilde{f}\left(Z_t\right) \big/ \partial \widetilde{Q}\left(t,Z_t\right)$ and $\widetilde{f}^{{\prime\prime}}\left(Z_t\right)= \partial^2 \widetilde{f}\left(Z_t\right) \big/ \partial \widetilde{Q}\left(t,Z_t\right)^2$, in the case where we have a time--dependent distribution function in our mapping, and we drop the dependence on $t$ in these functions otherwise.  Here, $Q(x)=F^{-}(x)$ and when referring to the ``true law'' of the driving process, we replace $Q$ with $Q_Y$, $F$ with $F_Y$ and $f$ with $f_Y$ to denote the quantile, distribution and density functions of this process, respectively.

\begin{Proposition}\label{gtransformlipschitzprop}
Let $(Z_t)$ be a $g$--transform quantile process given by Definition \ref{gtransformqp}, where the drift coefficient is given by either Eq. (\ref{gtransformdrift1}) or (\ref{gtransformdrift2}), and the volatility coefficient by Eq. (\ref{gtransformvol}). Let $(Y_t)_{0\le t<\infty}$ be a homogeneous driving process satisfying $\rd Y_t = \mu \rd t + \sigma \rd W_t$ for $\mu\in\mathbb{R}$, $\sigma\in\mathbb{R}^{+}$, and $Y_0=y_0\in\R$. Then the coefficients of $(Z_t)$ are Lipschitz continuous on $[-1/g,\infty)$, for $g\in(0,\infty)$, if the density function $f(t,y)$, associated with the law $F_{\bullet}(t,y)$, for all $t\in[t_0,\infty)$ has both left and right tail--decay to zero, is bounded on its support, and it is such that the following set of conditions is satisfied: 
\begin{align}\label{gliplimit1left}
        \lim_{z\rightarrow-{1/g}^+}\widetilde{f}\left(z\right)\left[g+\frac{\lnn(gz+1)}{g} \right]\expp\left(\frac{\left(\lnn(gz+1)\right)^2}{2g^2} \right) &= L_1<\infty, \\[0pt]\label{gliplimit1right}
 \lim_{z\rightarrow\infty}\widetilde{f}\left(z\right)\left[g+\frac{\lnn(gz+1)}{g} \right]\expp\left(\frac{\left(\lnn(gz+1)\right)^2}{2g^2} \right) &= L_2<\infty, \\[0pt]\label{gliplimit2left}
 \lim_{z\rightarrow-{1/g}^+} \widetilde{f}^{\prime}\left(z\right)\big/\widetilde{f}\left(z\right) &=  L_3<\infty, \\[0pt]\label{gliplimit2right}
 \lim_{z\rightarrow\infty}  \widetilde{f}^{\prime}\left(z\right)\big/\widetilde{f}\left(z\right) &= L_4<\infty, \\[0pt] \label{gliplimit3left}
\begin{split} \lim_{z\rightarrow-{1/g}^+}\widetilde{f}\left(z\right)^2\left[2\frac{\lnn^2(gz+1)}{g^2} + 3\lnn(gz+1) +(g^2+1) \right]\expp\left(\frac{\left(\lnn(gz+1) \right)^2}{g^2} \right) &= L_5<\infty,\end{split} \\[0pt] \label{gliplimit3right}
 \lim_{z\rightarrow\infty}\widetilde{f}\left(z\right)^2 \left[2\frac{\lnn^2(gz+1)}{g^2} + 3\lnn(gz+1) +(g^2+1) \right] \expp\left(\frac{\left(\lnn(gz+1) \right)^2}{g^2} \right) &=  L_6<\infty, 
\\[0pt]\label{gliplimit4left}
 \lim_{z\rightarrow-{1/g}^+}\widetilde{f}^{{\prime\prime}}\left(z\right)\big/\widetilde{f}\left(z\right) &= L_7<\infty, \\[0pt]\label{gliplimit4right}
  \lim_{z\rightarrow\infty}\widetilde{f}^{{\prime\prime}}\left(z\right)\big/\widetilde{f}\left(z\right) &= L_8<\infty.
  \end{align}
In the case where $g\in(-\infty,0)$, the conditions that must be satisfied are the same as the above, however the limits are taken as $z\rightarrow-\infty$ and $z\rightarrow-{1/g}^-$ to show Lipschitz continuity of the drift and volatility coefficients on $(-\infty, -1/g]$.
\end{Proposition}
\begin{proof}
In order for the drift and volatility coefficients given by Eqs (\ref{gtransformdrift1}) to (\ref{gtransformvol}) to be Lipschitz continuous, their first derivative must be bounded on the range on which they are differentiable everywhere.  The first derivatives with respect to $z$ of these expressions are given by
\begin{equation}\begin{split}\label{alphaderivnotlawg}
\frac{\partial}{\partial z}\alpha(t,z) =& \frac{\mu \widetilde{f}^{\prime}\left(z\right)}{\widetilde{f}\left(z\right)} + 2\pi\mu \widetilde{f}\left(z\right)\left[g+\frac{\lnn(gz+1)}{g} \right]\expp\left(\frac{\left(\lnn(gz+1) \right)^2}{2g^2} \right)  \\[0pt] 
&+ \frac{1}{2}\sigma^2 \left\{\frac{\widetilde{f}^{{\prime}{\prime}}\left(z\right)}{\widetilde{f}\left(z\right)} + 6\pi \widetilde{f}^{\prime}\left(Z_t\right)\left[g+\frac{\lnn(gz+1)}{g} \right] \expp\left(\frac{\left(\lnn(gz+1) \right)^2}{2g^2} \right)\right. \\[0pt] 
&\left.+2\pi  \widetilde{f}\left(z\right)^2 \left[2\frac{\lnn^2(gz+1)}{g} + 3\lnn(gz+1) +\left(g^2+1\right) \right]\expp\left(\frac{\left(\lnn(gz+1) \right)^2}{g^2} \right)\right\},
\end{split}\end{equation}
\begin{equation}\begin{split}\label{alphaderivtruelawg}
\frac{\partial}{\partial z}\alpha(t,z) =& \frac{1}{2}\sigma^2 \left\{\frac{2\widetilde{f}''\left(Z_t\right)}{\widetilde{f}\left(z\right)} 8\pi \widetilde{f}'\left(z\right)\left[g+\frac{\lnn(gz+1)}{g} \right]\expp\left(\frac{\left(\lnn(gz+1) \right)^2}{2g^2} \right) \right. \\[0pt] &\left.+2\pi  \widetilde{f}\left(z\right)^2 \left[2\frac{\lnn^2(gz+1)}{g} + 3\lnn(gz+1) +\left(g^2+1\right) \right]\expp\left(\frac{\left(\lnn(gz+1) \right)^2}{g^2} \right)\right\},
\end{split}\end{equation}
and 
\begin{equation}\begin{split}\label{sigmatildederivg}
\frac{\partial}{\partial z}\widetilde{\sigma}(t,z) =& \frac{\sigma \widetilde{f}^{\prime}\left(z\right)}{\widetilde{f}\left(z\right)} + 2\pi\sigma \widetilde{f}\left(z\right)\left[g+\frac{\lnn(gz+1)}{g} \right] \expp\left(\frac{\left(\lnn(gz+1)\right)^2}{2g^2} \right),
\end{split}\end{equation}
respectively.  

It follows that satisfying the limits in the Proposition ensures Eqs (\ref{alphaderivnotlawg}) to (\ref{sigmatildederivg}) are bounded for all $t\in[t_0,\infty)$ and $z\in[-1/g,\infty)$ for $g\in(0,\infty)$, or $z\in(-\infty,-1/g]$ for $g\in(-\infty,0)$. 
\end{proof}

\subsubsection{Skewed GBM and OU quantile diffusions}
We motivate the construction of a GBM and an OU $g$--transform quantile diffusion by the instance where we are interested in the quantiles of some process that has increased skewness relative to a process with a lognormal or a normal finite dimensional distribution, respectively.

First, consider a geometric Brownian motion (GBM) driving process satisfying the SDE $\rd Y_t = \mu Y_t\rd t + \sigma Y_t \rd W_t$ with $y_0\in\mathbb{R}^{+}$, $\mu\in\mathbb{R}$ and $\sigma\in\mathbb{R}^+$.  The ``true law'' GBM--Tukey--$g$ quantile diffusion, without loss of generality we will set the location and scale as standardised through $A=0$ and $B=1$ for simplicity, which then satisfies the SDE 
\begin{equation}\label{gbmgtransformsde}
\rd Z_t = \left( \dfrac{g}{2t} - \dfrac{\lnn(gZ_t+1)}{2gt}\right)\left(gZ_t+1\right) \rd t + \dfrac{\left(gZ_t+1\right)}{\sqrt{t}}\rd W_t
\end{equation}
for $t\in[t_0,\infty)$ with $Z_{t_0}=z_{t_0}\in[-1/g,\infty)$ when $g>0$ and $z_{t_0}\in(-\infty,-1/g]$ when $g<0$.  We emphasise that the $g$-transform quantile diffusion satisfying the SDE \ref{gbmgtransformsde} may also be constructed by taking $Y_t=W_t$ and transforming using the ``true law''---see Remark \ref{qpuniquerem}.

Now consider an Ornstein-Uhlenbeck (OU) driving process satisfying the SDE \newline  $\rd Y_t = \theta\left( \mu -Y_t\right)\rd t + \sigma  \rd W_t$ with $y_0\in\mathbb{R}^{+}$, $\mu\in\mathbb{R}$ and $\theta,\sigma\in\mathbb{R}^+$.  The ``true law'' OU--Tukey--$g$ quantile diffusion satisfies the SDE
\begin{equation}\label{ougtransformsde}
\rd Z_t = \left(\dfrac{g\theta}{\left(1-\expp\left(-2\theta t\right)\right)} - \dfrac{\theta\lnn(gZ_t+1)}{g\left(1-\expp\left(-2\theta t\right)\right)} \right)\left(gZ_t+1\right) \rd t + \dfrac{\sqrt{2\theta}\left(gZ_t+1\right)}{\sqrt{1-\expp\left(-2\theta t\right)}}\rd W_t
\end{equation}
for $t\in[t_0,\infty)$ with $Z_{t_0}=z_{t_0}\in[-1/g,\infty)$ when $g>0$ and $z_{t_0}\in(-\infty,-1/g]$ when $g<0$.

We emphasize that both SDEs, (\ref{gbmgtransformsde}) and (\ref{ougtransformsde}), can be written in the form
\begin{equation}
\rd Z_t = \dfrac{\sigma^2}{2{\rm{Var}}(Y_t\vert Y_0=y_0)}  \left(g - \frac{\lnn(gZ_t+1)}{g} \right)\left(gZ_t+1\right)\rd t + \dfrac{\sigma}{\sqrt{{\rm{Var}}(Y_t\vert Y_0=y_0)}}(gZ_t+1)\rd W_t, 
\end{equation}
where for any $t_0\leq t<\infty$ the diffusion coefficients are globally Lipschitz continuous on $z\in[-1/g,\infty)$ for $g>0$ and on $z\in(-\infty,-1/g]$ for $g<0$.  The first derivative of the drift coefficients of the SDEs are unbounded as $z\rightarrow\infty$ or $z\rightarrow-\infty$ in theses cases, respectively, and so there does not exist strong global solutions to the SDEs.

\subsection{Tukey--$h$ transform quantile diffusions}
We now introduce the $h$--transform class of distributions, where again we generate random variables with this distribution through a transformation the random variable $X$, however now allowing for flexible modelling of the heaviness of the tails through the parameter $h$. Denote the Tukey--$h$ transformation map by $T_{A,B,h}(x)$. The random variable $\phi_h$, the distribution of which will belong to the $h$--transform class, is defined by
\begin{equation}
\phi_h := T_{A,B,0,h}(X) := A + BX\expp\left(\frac{hX^2}{2}\right),
\end{equation}  
where $A\in\mathbb{R}, B\in\mathbb{R}^+$ are the location and scale parameters, respectively, and $h>0$ allows for flexible relative kurtosis of the transformed random variable. 

\begin{Definition}\label{htransformquantilefunction}
 The quantile function of the single parameter $h$--transform family of distributions is given by $Q_{\phi_h}\left(u;\bm{\xi}\right) = T_{A,B,0,h}(x_u)$ for quantile level $u\in[0,1]$ where $x_u=\sqrt{2}\erf^{-}(2u-1)$ is the $u^{th}$ quantile of the standard normal distribution and $\bm{\xi}=(A,B,h)$ is a vector of parameters.
\end{Definition}

\begin{figure}[H]
\setlength\fwidth{.825\textwidth}
\centering
\input{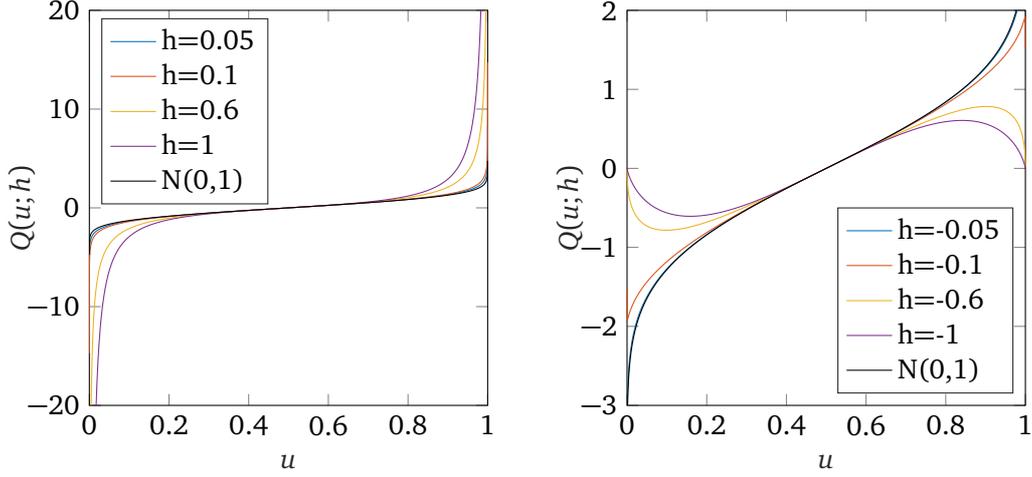}
\caption{$h$--transform quantile functions relative to the standard normal quantile function for $h\in\{-0.05, -0.1, 0.05, 0.1, 0.6, 1\}$, and $h$-quantile transforms for  $h\in\{-0.6, -1\}$.}
\label{htransformquantile}
\end{figure}

Figure \ref{htransformquantile} shows the $h$--transform quantile functions for varying values of the parameter $h$, relative to the quantile function $x_u$ of a standard normal random variable.  The plot on the right shows that for negative values of $h$ beyond a certain threshold, the quantile function is no longer monotonically increasing and hence we restrict to $h\in(0,\infty)$, i.e., we introduce more kurtosis to the base random variable.  In general it is assumed that $\vert h\vert< \approx 0.1$ so to prevent the tails becoming too pronounced and leading to almost instant growth.  For $h\in(0,\infty)$, we have $Q_{\phi_h}(u;A,B,h): [0,1] \rightarrow (-\infty,\infty)$.

\begin{Definition}\label{htransformqpdefn}
The random--level $h$--transform quantile diffusion is given by Definition \ref{processdrivendefn} where we take $Q_{\zeta}(u;\bm{\xi})$ to be $Q_{\phi_h}(u;A,B,h)$.
\end{Definition}

The dynamics of the $h$--transform quantile diffusion are derived in Appendix \ref{htransformcoefficients}. The cases under which the drift and diffusions coefficients derived here will ensure that a unique solution to the SDE \ref{dz1} exists can be characterised in terms of a set of conditions on the distribution function $F$ in the composite map, as stated in Appendix \ref{htransformlipconditions}.  

\subsection{Tukey $g$--$h$ quantile diffusions}
Finally, we consider the more general Tukey $g$--$h$ class of distributions, where again a random variable $\phi_{g,h}$ belonging to this family is generated through a transformation of the random variable $X$, given by
\begin{equation}\label{ghtransform}
\phi_{g,h}:=T_{A,B,g,h}(X)= A + B\left(\frac{\expp (gX) - 1}{g}\right)\expp\left(\frac{hX^2}{2}\right)
\end{equation}
where the parameters $g\in\mathbb{R}\setminus 0$ and $h\in\mathbb{R}^+$ are responsible for controlling the skewness and kurtosis of the distribution, respectively, and $A\in\mathbb{R}$ and $B\in\mathbb{R}^+$ are the location and scale parameters, respectively.  
\begin{Definition}\label{ghquantilefunction}
The quantile function of the $g$--$h$ family of distributions is given by $Q_{\phi_{gh}}(u;\bm{\xi}) = T_{A,B,g,h}(x_{u})$ for $u\in[0,1]$ where $x_u=\sqrt{2}\erf^{-}(2u-1)$ is the $u^{th}$ quantile of the standard normal distribution and $\bm{\xi}=(A,B,g,h)$ is a vector of parameters.
\end{Definition}
In either of the cases $g<0$ and $g>0$, we have $Q_{\phi_{gh}}\left(u;A,B,g,h\right):[0,1]\rightarrow (-\infty, \infty),$ where $h\geq0$ to ensure the quantile function is monotonically increasing.  
\begin{Definition}\label{ghprocessdriven}
The random--level $g$--$h$ quantile diffusion is given by Definition \ref{processdrivendefn} where we take $Q_{\zeta}(u;\bm{\xi})$ to be $Q_{\phi_{gh}}(u;A,B,g,h)$.
\end{Definition}
The drift and volatility functions of the $g$--$h$ quantile diffusion can be derived using Proposition \ref{prop1}.

\section{Application of quantile diffusions in distortion pricing}\label{examplessection}

In this section we consider an application that follows naturally from the properties of the random--level quantile diffusion.  The construction gives rise to an output (quantile) process whereby one has control over the higher--order moments induced by the transformation through direct parameterisation of such features, e.g., skewness and kurtosis, in the composite map.  Motivated in part by the argument made by Venter \cite{venter} that no-arbitrage pricing assumptions imply a distribution transformation, various distortion operators and transforms have been developed for the pricing of assets and risks, i.e., the distortions occur in such a way that the mean value under the distorted distribution yields a risk--adjusted premium or asset price.  We consider the Wang and Esscher transforms, both of which correspond to pricing in line with B{\"{u}}hlmann's economic principle \cite{buhlmann1980economic}, and a class of distortion operators given in \cite{godin2012} and generalised in \cite{godin2019}.  We present our approach based on quantile processes as a comparative and novel way to induce distribution distortions for pricing. 

In what follows, we allow the argument in the superscript of a distribution or quantile function to denote the probability measure under which the random variable has such distribution.

Let $(\Omega,\mathscr{F},\mathbb{P})$ be a probability space and $X$ be a random variable representing a financial risk with distribution function $F_{X}^{\mathbb{P}}(x)=\mathbb{P}(X\leq x)$.  Define the survival function of $X$ to be the strictly decreasing function $\bar{F}^\mathbb{P}_{X}=\mathbb{P}(X>x)$.   A distortion operator is an increasing, differentiable function $\nu:[0,1]\rightarrow[0,1]$ with $\nu(0)=0$ and $\nu(1)=1$, that acts on the CDF or survival function of a random variable.  We assume the existence of a measure $\mathbb{P}^{*}\sim\mathbb{P}$ on $(\Omega,\mathscr{F})$ such that the distorted $\mathbb{P}$--distribution of any random variable corresponds to its $\mathbb{P}^{*}$--distribution.

Among some of the most well--known distortion operators are the proportional hazard (PH) distortion function, see \cite{wang1996,wang1997axiomatic}, and the Wang transform \cite{wang,wang2}, which we present here.  Further pricing frameworks constructed from distortion operators are given in \cite{godin2012, godin2019,kijima,kijimamuromachi, kijimamuromachi2} and a summary of such methods is provided in \cite{cruz2015fundamental} and the references therein. 

Let $\Phi(x)$ be the standard normal CDF and $T_k(x)$ is the CDF of the Student--t distribution with location parameter $\mu=0$ and $k$ degrees--of--freedom.  The one--factor Wang distortion operator is given by $\nu_{\lambda}^{(1)}(u)= \Phi(\Phi^{-}(u)+\lambda)$ and the two--factor Wang distortion operator by $\nu_{\lambda}^{(2)}(u)=T_k(\Phi^{-}(u)+\lambda)$ where $u\in[0,1]$ and $\lambda\in\mathbb{R}$ is a parameter that captures the level of systematic risk.  Whilst the one--factor Wang transform is motivated from the perspective of its ability to recover CAPM and the Black--Scholes model, it has received criticism for its inability to produce a distorted measure that accounts for higher--order moments, e.g., heavy--tailed features or skewness, that are often observed in financial returns data.  The two--factor Wang transform \cite{wang2} overcomes such limitations, as well as accounting for parameter uncertainty that may arise when pricing under such a model.  Such a model is no longer consistent with the risk--neutral CAPM model or B\"{u}hlmann's pricing principle.  A further extension of the Wang transform---the generalised Wang transform---as developed in \cite{kijimamuromachi}, that is consistent with B\"{u}hlmann's principle and also provides larger flexibility in incorporating higher--order moments is given as follows.  

\begin{figure}[H]
\setlength\fwidth{\textwidth}
\centering
\resizebox{0.5\textwidth}{7.5cm}{\input{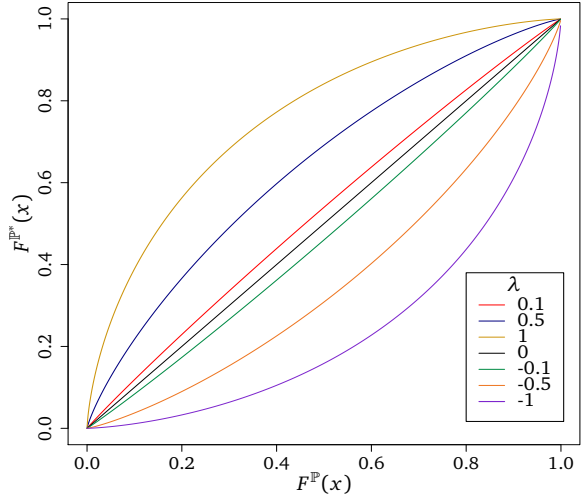}}
\caption{The one--factor Wang transform $F^{\mathbb{P}^{*}}(x)=\Phi(\Phi^{-}(F^{\mathbb{P}}(x))+\lambda)$ for a range of $\lambda$ parameters.}
\label{wangdistortionplot}
\end{figure}

Let $V$ be any positive valued random variable, $\phi$ a random variable with standard normal distribution, independent of $V$ and $G(x)$ the distribution of the random variable corresponding to the ratio $\phi/V$.  The generalised Wang distortion operator is given by $\nu_{\lambda}^{\rm{Gen}}=\mathbb{E}_V[\Phi(G^{-}(u) )V + \lambda ]$ where $u\in[0,1]$ and $\lambda\in\mathbb{R}$ represents the market price of risk.  In the case that the random variable $V$ is almost surely a Dirac mass on the event $\{V=1\}$, one recovers the standard Wang transform.  The risk--adjusted CDF of the random variable $X$ under the Wang transforms is given by \begin{equation}
F_X^{\mathbb{P}^{*}}(x):=\nu^{(i)}(F_X^{\mathbb{P}}(x))
\end{equation} 
where $i\in\{1, 2,\rm{Gen} \}$. The risk--adjusted actuarial fair--value of $X$, denoted $\pi(X)$, can be computed by taking the expectation of the discounted value $X$ with respect to $\PR^*$, that is $\pi (X)=\mathbb{E}^{\PR^*}[X]$.  Figure \ref{wangdistortionplot} illustrates the effect of the Wang transform on an input distribution $F^{\mathbb{P}}(x)$ for a range of $\lambda\in\R$ parameter values. 

The generalised class of distortion operators given in \cite{godin2019} assumes the existence of a probability measure $\mathbb{Q}$ on $(\Omega,\mathscr{F})$ that is equivalent to $\mathbb{P}$, such that under the conditions on $X$ given in \cite{godin2019}, the general distortion operator implied by the random variable $X$ is defined by 
\begin{equation}\label{godindistortionoperator}
g_X^{\mathbb{Q}\,\mathbb{P}}(u):=  \bar{F}_{X}^{\mathbb{Q}}\left( \bar{F}_{X}^{-\,\mathbb{P}}(u) \right) = 1-F_{X}^{\mathbb{Q}}\left(Q_{X}^{\mathbb{P}}(1-u) \right),
\end{equation}
for $u\in[0,1]$.  This class of distortion operators is flexible in its ability to incorporate higher--order distribution features such as skewness and kurtosis, and allows one to produce prices that are consistent with no-arbitrage models, equilibrium models and actuarial pricing principles. In the case that Eq. (\ref{godindistortionoperator}) is used to price financial derivative contracts (with continuous and increasing payoff functions) written on the random variable $X$, the distorted $\mathbb{P}$--distribution of the contract coincides with its $\mathbb{Q}$-distribution, where $\mathbb{Q}$ is a given pricing measure such as a risk--neutral measure. The distortion operator used is that implied by the financial risk that underlies the derivative contract, i.e., Eq. (\ref{godindistortionoperator}) when the value of the asset underlying the financial derivative is modelled by the random variable $X$. The connection with the Radon--Nikodym derivative is also given.   

Finally, the Esscher transform, which is based on weighting the marginal CDF of the financial risk $X$, is defined by
\begin{equation}\label{esscherdistributionfunction}
\rd F_X^{\mathbb{P}^{*}}(x):=\dfrac{\expp\left(\lambda x\right)\rd F_X(x)}{\int_{\mathbb{R}}\expp\left(\lambda y\right)\rd F_X(y)}
\end{equation}
where $\lambda\in\mathbb{R}$ is a risk--aversion parameter.

We note that the above pricing distortions do not rely on the assumption of market completeness.

\subsection{Distorted measures induced by stochastic quantile transformations}

The construction of quantile diffusions presented in this paper is motivated largely from the perspective of the induced statistical properties of the class of models, and the direct control one has over such features through the composite map.  We consider an application to distortion--based pricing that utilises the induced measure of the finite dimensional distributions of quantile diffusions.  The functional form of the quantile process transformation allows the models to encompass the induced higher--order features as an easily interpretable, and general, parametric family.  The models allow for the ability to change the moments non--trivially in a targeted way so to encompass risk preferences and market risk exposure that an investor much be compensated for, much like the parameters in the above distortion operators.  An advantage lies in the unrestricted nature of the class of parametric quantile models that may be developed, and that the interplay between the statistical properties of the driving process and the output quantile process is often explicit.  We note that whilst parallels can be drawn between the properties of each approach, the composite map used in the quantile diffusion construction does not define a distortion operator.  Additionally, the following distorted measure framework allows for the dynamic pricing of risks in continuous time.

We introduce the measure distortion based on quantile diffusions as follows. In what follows, we omit the dependence on parameters in the notation for distribution, quantile and density functions.

Consider a probability space $(\Omega,\mathscr{F},\mathbb{P})$ equipped with filtration $(\mathscr{F}_t)_{0\leq t<\infty}$.  We have $(Y_t)_{0\leq t<\infty}$ a diffusion process on the probability space with marginal distribution function $F_Y^{\mathbb{P}}(t,y):\mathbb{R}^{+}\times D_Y\rightarrow[0,1]$ for $t>0$ and $D_Y:=[{\rm ess \hspace{1mm} inf}\hspace{1mm} Y_t, {\rm ess \hspace{1mm} sup}\hspace{1mm} Y_t]\subseteq\mathbb{R}$, $Y_0=y_0\in D_Y$, and assume that $\mathscr{F}_t:=\sigma((Y_s)_{0\leq s\leq t})$.  Denote the marginal quantile function of the process at each $t\in(0,\infty)$ by $Q_Y^{\mathbb{P}}(t,y):=F_Y^{\mathbb{P}-}(t,y)$ for $u\in[0,1]$. We write the transition distribution function of $(Y_t)$ as
\begin{equation}\label{yconditionaldist}
F_Y^{\mathbb{P}}(t,y\vert\, s,y_s) = \mathbb{P}\left(Y_t\leq y\vert\, Y_s=y_s\right)
\end{equation}
for $y_s, y\in D_Y$ and for all $0<s\leq t<\infty$. When $s=0$ and $t>0$, this reduces to the marginal case.
\begin{Definition}
Consider Definition \ref{processdrivendefn} where $Q_{\zeta}(u)$ is the quantile function of some random variable $\zeta:\Omega\rightarrow D_{\zeta}\subseteq\mathbb{R}$ and $(Y_t)$ is the driving process.  We construct a random--level quantile diffusion $(Z_t)_{t_0\leq t<\infty}$ as in case (i) of the definition.  The process has marginal distribution function under the $\mathbb{P}$-measure given by $F_Z^{\mathbb{P}}(t,z)= \mathbb{P}\left(Z_t\leq z\right) = \mathbb{P}\left(Q_{\zeta}(F_Y^{\mathbb{P}}(t,Y_t))\leq z\right)$, for $z\in D_{\zeta}$ and $t\in[t_0,\infty)$.  We may also write the transition distribution function of $(Z_t)$ as
\begin{equation}\begin{split}\label{zconditionaldist}
 F_Z^{\mathbb{P}}(t,z\vert\, s, z_s) = \mathbb{P}\left(Z_t\leq z\vert Z_s = z_s \right) = \mathbb{P}\left(Q_{\zeta}\left(F_Y^{\mathbb{P}}\left(t,Y_t\right) \right)\leq z\vert Y_s=y_s\right)
\end{split}\end{equation}
for all $z\in D_{\zeta}$, where $z_s=Q_{\zeta}(F_Y^{\mathbb{P}}(s,y_s))\in D_{\zeta}$ and for all $t_0\leq s\leq t<\infty$.
\end{Definition}
If the composite map relating $(Y_t)$ and $(Z_t)$ is invertible and the transition distribution of $(Y_t)$ is known, that of $(Z_t)$ can be expressed as $F_Z^{\mathbb{P}}(t,z\vert\, s,z_s)= F_Y^{\mathbb{P}}(t,Q_Y^{\mathbb{P}}(t,F_{\zeta}(z)) \vert\, s, Q_Y^{\mathbb{P}}(s,F_{\zeta}(z_s)))$, given that probability mass is conserved.
\begin{rem} The focus here is on using the quantile transformation of $(Y_t)$ to construct a new, equivalent measure $\mathbb{P}^{Z}$ on the space $(\Omega, \mathscr{F})$, whereby the quantile transformation reassigns the probabilities in $\Omega$ under $\mathbb{P}$ to those under $\mathbb{P}^{Z}$.  Such a quantile measure--transform defines a Radon-Nikodym derivative.
\end{rem}

\begin{Definition}\label{distorteddynamicmeasure}
Let $(Y_t)$ and $(Z_t)$ be the driving and quantile diffusions, respectively, and $D_Y=D_{\zeta}\subseteq\mathbb{R}$.  We assume the existence of a new measure $\mathbb{P}^{Z}\sim\mathbb{P}$ such that for $t_0\leq s\leq t<\infty$, the transition distribution of the driving process under $\mathbb{P}^{Z}$ is given by
\begin{equation}
F_Y^{\mathbb{P}^{Z}}(t,y\vert\, s,y_s)=F_Z^{\mathbb{P}}(t,y\vert\, s,y_s),
\end{equation}
 for $y\in D_Y=D_{\zeta}$, where $Z_s=y_s\in D_{\zeta}$ is the value of the quantile diffusion at time $s$ that we observe under $\mathbb{P}$. Marginally, we define $F_Y^{\mathbb{P}^{Z}}(t,y)=F_Z^{\mathbb{P}}(t,y)$ for $y\in D_Y$ and $t\in[t_0,\infty)$. The conditional distorted measure is defined for each $t\in(s,T]$ for $0<t_0\leq s<T<\infty$ by
\begin{equation}\label{conditionaldistortedmeasure}
\mathbb{P}_{s,t}^{Z}(A\,\vert\,\mathscr{F}_s):= \int_{y\in A}\rd F_Y^{\mathbb{P}^{Z}}\left(t,y\vert\, s,y_s\right) = \int_{y\in A} \rd F_Z^{\mathbb{P}}\left(t,y\vert\, s,y_s\right)
\end{equation}
for all $A\subseteq D_Y$. Unconditionally, the distorted measure is defined for each $t\in[t_0,T]$, by
\begin{equation}
    \mathbb{P}_t^Z(A):= \int_{y\in A}\rd F_Y^{\mathbb{P}^Z}(t,y) = \int_{y\in A}\rd F_Z^{\mathbb{P}}(t,y)
\end{equation}
for all $A\subseteq D_Y$.  When $t=T$, we write $\mathbb{P}_T^Z(A)=\mathbb{P}^Z(A)$ for all $A\subseteq D_Y$.  It follows that for $t\in(0,\infty)$ and $Y_t$ an $\mathscr{F}_t$--adapted random variable, we have
\begin{enumerate}[(i)]
    \item $\mathbb{E}^{\mathbb{P}^Z}[Y_t] := \mathbb{E}^{\mathbb{P}_t^Z}[Y_t]\;$,
    \item $\mathbb{E}^{\mathbb{P}_{s,t}^Z}\left[Y_t\right]:=\mathbb{E}^{\mathbb{P}_t^{Z}}\left[Y_t\vert\mathscr{F}_s\right],\;$ for $\; 0<s<t\leq T<\infty$.
\end{enumerate}
\end{Definition}
It follows that for $\varphi(y)$ any Borel function and $t\in(s,T]$,
\begin{equation}\label{conditionalexpmeasurechange}
\mathbb{E}^{\mathbb{P}_{s,t}^Z}[\varphi(Y_t)]=\mathbb{E}^{\mathbb{P}_t^{Z}}\left[\varphi(Y_t)\vert\mathscr{F}_s \right] = \mathbb{E}^{\mathbb{P}}\left[\varphi(Z_t)\vert\mathscr{F}_s \right] = \mathbb{E}^{\mathbb{P}}\left[\varphi\left(Q_{\zeta}\left(F_Y^{\mathbb{P}}\left(t,Y_t\right) \right)\right)\vert\mathscr{F}_s \right]. 
\end{equation}
The above equation is useful when perhaps we want to value a contingent claim written on $(Y_t)$ with price process $\varphi_t=\varphi(Y_t)$ under a probability measure that accounts for risk--inducing skewness and kurtosis.  Here, $\mathbb{P}^{Z}$ might be the measure induced by the distribution of a Tukey $g$--$h$ quantile diffusion $(Z_t)$ that is constructed by considering $(Y_t)$ to be the driving process, for example.   This quantile diffusion characterises the law of an implicit process with the desired levels of skewness and kurtosis we wish to capture in the valuation measure.  Eq. (\ref{conditionalexpmeasurechange}) would then give the (non-discounted) value of the contingent claim at time $s\in[t_0,T]$.

The existence of $\mathbb{P}^{Z}$ follows from the fact that each distribution function induces a unique probability measure on the Borel $\sigma$-algebra of $\mathbb{R}$.  The mechanism that takes us from $\mathbb{P}$ to $\mathbb{P}^{Z}$ is the quantile transformation.  In other words, we view the probabilities associated to the distorted random variable $Z_t$ for each $\omega\in\Omega$, $t\in[t_0,T]$ under $\mathbb{P}$ as those assigned to the random variable $Y_t$ for each $\omega\in\Omega$, $t\in[t_0,T]$, under $\mathbb{P}_{t}^{Z}$.  The $\mathbb{P}^{Z}$ measure redistributes the probabilities in $\Omega$ to account for the properties, e.g., the skewness and kurtosis, we factor in to the quantile transformation.  The assumption that $D_Y=D_{\zeta}$ is necessary to ensure the probability measures are equivalent.  When applied to the valuation of a financial risk, the distortion aims to capture the risk exposure that must be compensated for.  

\begin{rem}
Interpreting the distributional distortion induced by quantile diffusions as that inducing a measure change is equivalent to constructing a measure change with the focus of inducing particular statistical properties under the new measure.  The dynamics of $(Y_t)$ under the new measure will be equivalent to those of $(Z_t)$ under the original measure. 
\end{rem}
\begin{Example}
The purpose of this example is to demonstrate that the well established Wang distortion transforms, widely considered in actuarial science, may be replicated as special cases of the general distortion framework proposed in our work. One noteworthy aspect of this is that the Wang transform approach is treated in a non--dynamic manner in the literature. As such, one can also see this following example as both capturing the general family of Wang distortion measures as well as generalising them into a consistent dynamic framework.

Consider Definition \ref{processdrivendefn} (ii) and construct a quantile diffusion using ``false law'' $F(t,y;\lambda)=\Phi(y+\lambda)$ for $\lambda\in\mathbb{R}$, and stationary driving process $(Y_t)$ that is marginally standard normally distributed, that is $F_Y(t,y)=\Phi(y)$, for all $t\in(0,\infty)$.  For each $t\in(0,\infty)$, we have $Z_t:=Q_{\zeta}(\Phi(Y_t+\lambda);\bm{\xi})$ and therefore one can express the distribution as follows: $F_Z^{\mathbb{P}}(t,z;\bm{\xi},\lambda)=\mathbb{P}(Z_t\leq z)=\mathbb{P}(Q_{\zeta}(\Phi(Y_t+\lambda);\bm{\xi})\leq z)=\mathbb{P}(Y_t\leq \Phi^{-}(F_{\zeta}(z;\bm{\xi}))+\lambda)=\Phi(\Phi^{-}(F_{\zeta}(z;\bm{\xi}))+\lambda))$, for $z\in D_{\zeta}$.  Therefore, the marginal distorted distribution, i.e., that induced by the quantile process, is equivalent to the distortion induced by a one--factor Wang transform acting on the base distribution function $F_{\zeta}=Q_{\zeta}^{-}$.  The two factor Wang transform is replicated similarly if we consider a driving process with stationary marginal distribution  $F_Y(t,y;\bm{\theta})=T_k(y)$ for all $t\in(0,\infty)$.  
\end{Example}
We derive the Radon--Nikodym derivative process for the change of measure defined in Definition \ref{distorteddynamicmeasure} as follows.  We note that the measure induced by the distribution of the quantile diffusion cannot be achieved by a Girsanov transformation as the volatility coefficient, given in Eq. (\ref{volatility}) of the SDE satisfied by the quantile process differs to that of the driving process under all transformations that are not the identity map.
\begin{Definition}
Consider the filtered probability space $(\Omega,\mathscr{F},(\mathscr{F}_t),\mathbb{P})$, and the measure $\mathbb{P}^Z$ on the space such that it is the distorted measure induced by the marginal law of $Z_T$, as given in Definition \ref{distorteddynamicmeasure}.  By the Radon--Nikodym theorem, since $\mathbb{P}^{Z}\sim\mathbb{P}$ is assumed in the construction of $\mathbb{P}^Z$, there exists a non-negative random variable $\rho(\omega):\Omega\rightarrow\mathbb{R}_0^+$ such that the distorted measure at time $T$ is given by 
\begin{equation}\label{rndistortedmeasure}
\mathbb{P}^{Z}(A):= \int_{\{ \omega\in\Omega\hspace{1mm}:\hspace{1mm}Z_T(\omega)\in A\}} \rho (\omega) \rd \mathbb{P}(\omega)
\end{equation}
for all $A\in\mathscr{F}$.  This may be written as $\rd\mathbb{P}(\omega)\rho(\omega)=\rd\mathbb{P}^{Z}(\omega)$ for all $\omega\in\Omega$, and $\rho$ is called the Radon--Nikodym derivative of $\mathbb{P}^{Z}$ with respect to $\mathbb{P}$.  
Let $\rho$ be an almost--surely positive random variable satisfying $\mathbb{E}^{\mathbb{P}}[\rho]=1$. The Radon--Nikodym derivative process is defined by
\begin{equation}
\rho_t := \mathbb{E}^{\mathbb{P}}\left[\rho \vert \mathscr{F}_t \right] = \dfrac{\rd\mathbb{P}^{Z}}{\rd \mathbb{P}}\Big\vert_{\mathscr{F}_t}
\end{equation}
for all $t\in[t_0,T]$.
\end{Definition}
Moreover, assuming the existence of transition densities for the processes at each $t_0\leq s\leq t\leq T<\infty$ under both the $\mathbb{P}$ and $\mathbb{P}_{ s,t}^{Z}$ measures, and given the transition distribution functions in Eqs. (\ref{yconditionaldist}) and (\ref{zconditionaldist}), the conditional likelihood ratio, conditional on $(\mathscr{F}_s)_{t_0\leq s \leq t\leq T<\infty}$ is given by
\begin{equation}
\psi_Y^{\mathbb{P},\mathbb{P}^{Z}}(t,y \vert\, s,y_s):= \dfrac{\partial F_Y^{\mathbb{P}^{Z}}(t,y\vert\, s,y_s)}{
\partial F_Y^{\mathbb{P}}(t,y\vert\, s,y_s )}.
\end{equation}
for all $y\in D_Y$ and where $ y_s\in D_Y$, $z_s=Q_{\zeta}(F_Y^{\mathbb{P}}(s,y_s))$.

To summarise, we induce measure distortions by constructing a quantile diffusion from some driving diffusion, and then deriving the measure induced by the probability distribution of this quantile process.  The focus is to allow the distortion to capture the risk exposure that the investor must be compensated for, when valuing under the new measure, in such a way that this risk is parameterised directly in the construction of the model.  This section provides a framework that constructs such a measure $\mathbb{P}^{Z}$ capturing these higher--order risks through the construction of a random--level quantile diffusion.  Taking expectations of a Borel function of the driving process under $\mathbb{P}^{Z}$ will be equivalent to taking expectations of the output quantile process under the original measure.
\subsubsection{Measure--change under a Tukey $g$--$h$ quantile transformation}
We next consider the driving process to be a univariate drifted Brownian motion under $\mathbb{P}$ satisfying $\rd Y_t=\mu \rd t +\rd W_t$, $\mu\in\mathbb{R}$ with $Y_0=y_0\in\mathbb{R}$ for $t\in(0,T]$, $T<\infty$. 

\begin{figure}[H]
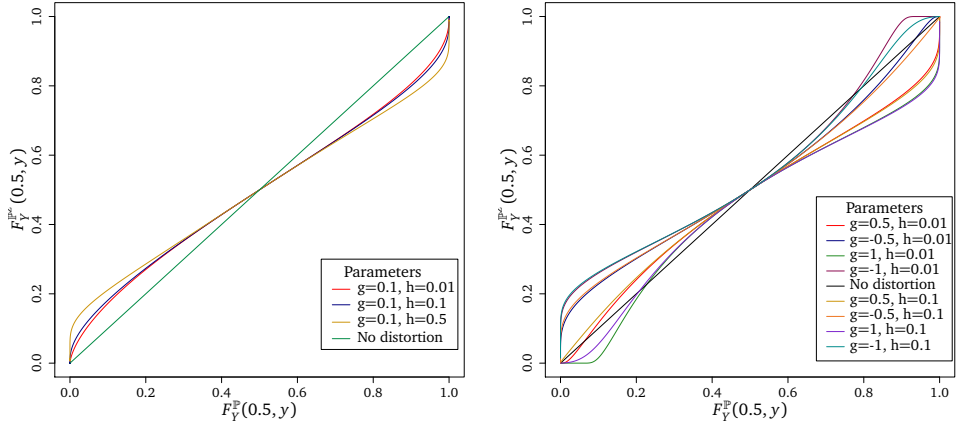

\centering
\resizebox{0.4\textwidth}{0.4\textwidth}{\input{distortionplotgh1a.tex}}
\resizebox{0.4\textwidth}{0.4\textwidth}{\input{distortionplotgh1b.tex}}
\vspace{-.45cm}
\caption{Relation between the marginal distribution functions of $(Y_t)$ at $t=0.5$ for $\mu=0$, under the $\mathbb{P}$ and $\mathbb{P}^{Z}$ measures, for a range of $g$ and $h$ parameters.}
\label{distortionplots1}
\end{figure}

Let $\zeta=\phi_{gh}$ and $Q_{ \phi_{gh}}(u;\bm{\xi})$ be the Tukey $g$--$h$ quantile function for $u\in[0,1]$, given by Eq. (\ref{ghquantilefunction}), where here $A=0$ and $B=1$ in the vector of parameters $\bm{\xi}=(A,B,g,h)$. Then, for $t\in[t_0,\infty)$, $t_0>0$, the random--level quantile diffusion is given by $Z_t=Q_{\phi_{gh}}(F_Y^{\mathbb{P}}(t,Y_t))$.  
The $\mathbb{P}$--transition law of $(Z_t)$ for $t_0\leq s\leq t<\infty$ is defined by Eq. (\ref{zconditionaldist}) as
\begin{equation}\begin{split}
F_Z^{\mathbb{P}}(t,z \vert\, s,z_s) 
&= \frac{1}{2}\left[1+\erf\left(\dfrac{\sqrt{2t}\hspace{1mm}\erf^{-}(2F_{\phi_{gh}}(z)-1)-\sqrt{2s}\hspace{1mm}\erf^{-}(2F_{\phi_{gh}}(z_s)-1)}{\sqrt{2(t-s)}} \right) \right]
\end{split}\end{equation}
where $F_{\phi_{gh}}(z):=Q^{-}_{\phi_{gh}}(z)$. By Definition \ref{distorteddynamicmeasure}, we set $F_Y^{\mathbb{P}^{Z}}(t,y \vert\, s,y_s)=F_Z^{\mathbb{P}}(t,y\vert\, s,y_s)$ for $y_s, y\in \mathbb{R}$, and define the Tukey $g$--$h$ conditional distorted measure by Eq. (\ref{conditionaldistortedmeasure}) for each $t\in(s,T]$, that is
\begin{equation}
\mathbb{P}^{Z}_{ s,t}\left(A\,\vert\,\mathscr{F}_s\right)=\int_{y\in A} \rd  \frac{1}{2}\left[1+\erf\left(\dfrac{\sqrt{2t}\hspace{1mm}\erf^{-}(2F_{\phi_{gh}}(y)-1)-\sqrt{2s}\hspace{1mm}\erf^{-}(2F_{\phi_{gh}}(y_s)-1)}{\sqrt{2(t-s)}} \right) \right],
\end{equation}
for all $A\subseteq \mathbb{R}$. The dependence on the skewness and kurtosis parameters, given by $g$ and $h$, respectively, comes in through the function $F_{\phi_{gh}}(y)$ which can be computed analytically.  The relation between the distribution functions that induce the $\mathbb{P}$-- and $\mathbb{P}^{Z}$--measures in this example are illustrated in Figures \ref{distortionplots1} and \ref{distortionplots2} for a range of the $g$ and $h$ parameters.

\begin{figure}[H]
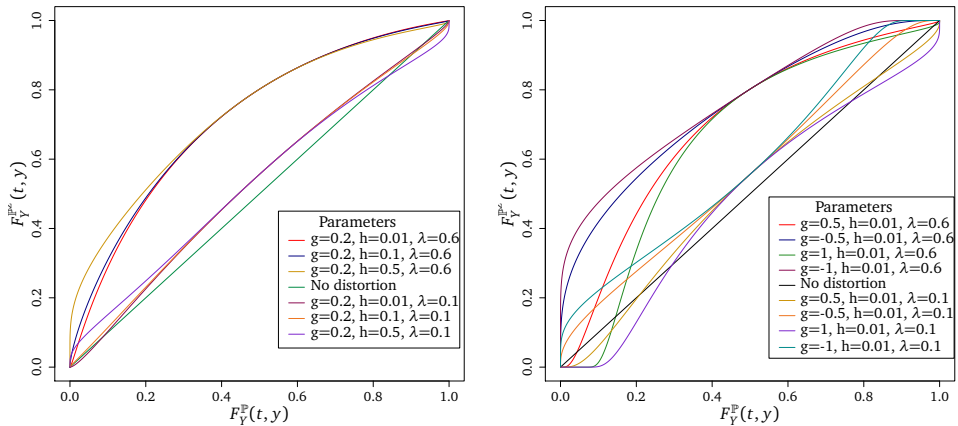

\centering
\resizebox{0.4\textwidth}{0.4\textwidth}{\input{distortionplotghfalse.tex}}
\resizebox{0.4\textwidth}{0.4\textwidth}{\input{distortionplotghfalse2.tex}}
\vspace{-.25cm}
\caption{Relation between the marginal distribution functions of $(Y_t)$ at $t=0.5$ for $\mu=0$, under the $\mathbb{P}$ and $\mathbb{P}^{Z}$ measures when we use a false law $F^{\mathbb{P}}(x)=\Phi(x+\lambda)$ in the mapping, for a range of $g$, $h$ and $\lambda$ parameters.}
\label{distortionplotsfalse}
\end{figure}

Assume we want to construct a Tukey $g$--$h$ quantile diffusion, however now we use the ``false law'' in the map---Definition \ref{processdrivendefn} (ii).  To draw a parallel with the Wang transform, let the ``false law'' $F^{\mathbb{P}}(x)$ be given by the time--homogeneous $\mathscr{N}(\lambda,1)$ normal distribution, for $\lambda\in\mathbb{R}$ and denote the quantile function of this distribution by $Q^{\mathbb{P}}(u)$ for $u\in[0,1]$.  It follows that for all $t\in[t_0,\infty)$, we construct the quantile diffusion by $Z_t=Q_{\phi_{gh}}(F^{\mathbb{P}}(Y_t))$ and derive its marginal distribution as $F_Z^{\mathbb{P}}(t,z) = \mathbb{P}(Z_t\leq z) = \mathbb{P}(Q_{\phi_{gh}}(F^{\mathbb{P}}(Y_t))\leq z ) = F_Y^{\mathbb{P}}(t, Q^{\mathbb{P}}(F_{\phi_{gh}}(z)))$ for $z\in\mathbb{R}$.  The distorted marginal law of the driving process is then given by $F_Y^{\mathbb{P}^{Z}}(t,y)=F_Z^{\mathbb{P}}(t,y) =F_Y^{\mathbb{P}}(t, Q^{\mathbb{P}}(F_{\phi_{gh}}(y)))$ for all $t\in[t_0,\infty)$ and $y\in\mathbb{R}$.  Figure \ref{distortionplotsfalse} illustrates the relation between the distorted and non--distorted marginal distributions for a range of $g$, $h$ and $\lambda$ parameters.
\vspace{-.5cm}
\vspace{-1cm}
\begin{figure}[H]
\centering
\resizebox{0.4\textwidth}{0.4\textwidth}{\input{distortionplotgh1c.tex}}
\resizebox{0.4\textwidth}{0.4\textwidth}{\input{distortionplotgh1d.tex}}
\vspace{-.5cm}
\caption{Relation between the marginal distribution functions of $(Y_t)$ at $t=0.5$ for $\mu=0.8$, under the $\mathbb{P}$ and $\mathbb{P}^{Z}$ measures, for a range of $g$ and $h$ parameters.}
\label{distortionplots2}
\end{figure}
\vspace{-1cm}

\subsection{Numerical example}
Consider Example 1 in \cite{wang}. Let a risk process $(Y_t)_{0\leq t<T}$ be such that at $T=1$ the random variable $Y_T$ has finite dimensional distribution given by the Pareto distribution function 
\begin{equation}
    F_Y(T,y)=F_Y(y)=1-\left(\dfrac{2000}{2000+y}\right)^{1.2},
\end{equation}
for $y>0$.  Consider the false--law shifted Tukey--$g$ quantile process with $(Y_t)$ as the driving process, given by $Z_t=B\expp(gY_t+g\gamma t)/g$ for all $t\in[t_0,T]$ and $B,g>0, \gamma\in\mathbb{R}$.  Let $\mathbb{P}^Z$ be the measure induced by the distribution of the random variable $Z_T$.  We compare the Tukey--$g$ risk--adjusted premium by layer to that produced under the Wang transform with parameter $\lambda=0.1$ and PH--transform with index $r=0.9245$.  Prices are obtained by taking expectations under $\mathbb{P}^Z$, that is $\pi_{T_g}[Y_T]=\mathbb{E}^{\mathbb{P}^Z}[Y_T]$, muchlike the pricing formula $H[X;\lambda]$ in \cite{wang}.    

We allow $g>0$ to introduce skewness relative to the risk distribution to capture investor risk--aversion levels, to account for skewness--related systemic risk and to place probabilistic emphasis on the occurrence of higher losses under the transformation.  The results are presented in Table \ref{gPCPtable}, where in the third column, the false--law location parameter $\gamma$ is selected to match the premium with the Wang and PH premia for the basic limit layer $(\$0,\$50000]$, as shown in bold, illustrating how the quantile transform influences the premiums for subsequent layers relative to this figure.  The same is done in the fourth column, however $\gamma$ is now chosen so to match the premiums for the layer $(\$200,\$300]$, shown in bold.  We observe that prior to this layer, the Tukey--$g$ premium is lower than the Wang and PH premia, and higher in subsequent layers with significantly larger layer--by--layer increases towards the higher layers resulting from the introduction of relative skewness under the quantile transform.  Overall, we observe that the PH--transform premium increases faster than the Wang transform, and the Tukey--$g$ measure distortion produces a premium that increases much faster than the PH--transform with the rate of increase determined by the magnitude of $g$.  

\begin{table}[H]
  \centering
 \begin{tabular}{|c|c|c|c|c|}
 \hline
  \multirow{2}{*}{Layer in 000's} & PH premium & Wang premium & \multicolumn{2}{c|}{$\pi_{T_g}[Y_{1};\gamma,g,B=0.01]$} \\
  \cline{4-5}
  & $r=0.9245$ & $\alpha=0.1$ & $g=0.8,\gamma=6.96$ & $g=0.08,\gamma=-10.25$ \\
\hline
(0,50] & \bf{5,487.0} & \bf{5,487.0} & \bf{5,487.2} & 261.9  
\\
\hline
(50,100] & 910 & 845.0 & 4,632.6 & 228.5  \\
\hline
(100,200] & 857 & 769.9 & 8,642.6 & 431.7  \\
\hline
(200,300] & 475 & \bf{414.2} & 8,219.7 & \bf{414.2} 
\\
\hline
(300,400] & 325 & 278.4 & 7,965.5 & 403.5  \\
\hline
(400,500] & 246 & 207.3 & 7,785.9 & 395.9  \\
\hline
(500,1000] & 728 & 598.0 & 37,287.4 & 1,909.3  \\
\hline
(1000,2000] & 675 & 533.2 & 70,349.5 & 3,635.2  \\
\hline
(2000,5000] & 819 & 616.6 & 197,478.0 & 10,306.4  \\
\hline
(5000,10000] & 567 & 405.7 & 310,290.9 & 16,331.72  \\
\hline
\end{tabular}
 \caption{Risk--adjusted premiums by layer under the Wang transform and Tukey--$g$ distorted measure.  The figures in bold correspond to the premiums matched across the different distortions.  The Tukey--$g$ premium for the layer $(0,50]$ has been matched (through  the choice of $\gamma$) to that of the PH and Wang premia for the same layer.  The Tukey--$g$ premium for the layer $(200,300]$ is matched to that of the Wang transform for the same layer.}
  \label{gPCPtable}
\end{table}

\section{Conclusions}

In this work we propose a novel, mathematical framework for the development of continuous--time, dynamic distortions induced by quantile processes. This approach can as such be widely adopted in dynamic risk analysis and applications. An important result of this framework is that one obtains a time--consistent method for the construction of interpretable, continuous--time measure distortion flows on the space of quantile functions. We demonstrate that the framework generalises, but includes as special cases, existing well--known measure distortion approaches used in financial mathematics, financial risk and actuarial science, while extending these approaches into dynamical settings. Furthermore, the developed framework is constructed from families of transformations, which produce statistically interpretable measure distortions induced by the quantile processes, making their application directly interpretable in relation to aspects of the risk inherent to the underlying driving process being accounted for in the distortion map. This has significant meaning in, among many other settings, financial risk management. For instance, it can be seen as a generalisation of classical risk pricing, often based on risk premium related to trend or volatility, that instead accommodates pricing in higher--order risk associated to skewness or kurtosis of the underlying risky process being assessed.

We developed two methods for the construction of the quantile processes to achieve our objectives: the first, which we largely focus on, features a dynamic random quantile level and allows for direct interpretation of the resulting quantile process characteristics such as location, scale, skewness and kurtosis, in terms of the model parameters. The second type are function--valued quantile diffusions and are driven by stochastic parameter processes, which determine the entire quantile function at each point in time. We derived core results relating to the definition, existence and uniqueness of such processes, and established the construction of the resulting distortion measure flows on the quantile space. We then produced key examples based on the flexible Tukey family of transforms applied to widely used underlying stochastic processes such as geometric Brownian motion and the Ornstein--Uhlenbeck processes. Such processes can be found in a range of practical application domains that include dynamic risk measures in econometrics, operations research sequential decision making, information theory and signal processing, and not least in general risk theory and applications thereof, e.g., in the context of dynamic risk assessment of climate change processes.

\vspace{.5cm}
\noindent{\bf Acknowledgments}. The authors acknowledge financial support by UK Research \& Innovation via EPSRC CASE project award 1939295, and are grateful for comments and suggestions by G. Kassis, W. T. Shaw, M. Crowe, and participants in the Mathematics \& Finance Conference, Research in Options, IMPA, Rio de Janeiro, Brazil, and Khalifa University, Abu Dhabi, UAE (December 2020), and attendees of the AIFMRM Seminar Series, University of Cape Town (September 2021). 
\vspace{-.25cm}
\bibliography{QD1refs}
\vspace{-.25cm}
\begin{appendix}
\section{Appendix}
\subsection{Existence of strong and weak solutions}\label{strongsolutionsection}
We consider the SDEs satisfied by diffusions, with a view to describe the conditions under which solutions to these SDEs exist, and whether they are strong or weak. 
Consider a probability space $(\Omega, \mathscr{F},\mathbb{P})$ and an $r$-dimensional Brownian motion along with its natural filtration,
\begin{equation}\label{rdimbm}
W=(\bm{W}_t, \mathscr{F}^W_t; 0\leq t<\infty ),
\end{equation}
where we assume the space is rich enough to accommodate a random vector $\bm{\xi}\in\mathbb{R}^d$, independent of $\mathscr{F}_\infty^W$ and with given distribution $\mu(\Gamma)=\mathbb{P}(\bm{\xi}\in\Gamma)$, $\Gamma\in\mathscr{B}(\mathbb{R}^d)$.  Moreover, we consider the left--continuous filtration $\mathscr{G}_t\triangleq\sigma(\bm{\xi})\vee\mathscr{F}_t^W = \sigma(\bm{\xi},W_s; 0\leq t<\infty)$ and the collection of null sets $\mathscr{N}\triangleq\{N\subseteq\Omega; \exists G\in\mathscr{G}_\infty\hspace{1mm}\textrm{with}\hspace{1mm}N\subseteq G\hspace{1mm}\textrm{and}\hspace{1mm}\mathbb{P}(G)=0 \}$ to obtain the augmented filtration 
\begin{equation}\label{augfiltration}
\mathscr{F}_t\triangleq\sigma\left(\mathscr{G}_t\cup\mathscr{N}\right), \hspace{1mm}0\leq t<\infty, \hspace{1mm} \mathscr{F}_\infty\triangleq\sigma\left(\cup_{t\geq0}\mathscr{F}_t\right).
\end{equation}
Next, we consider the general $d$--dimensional SDE
\begin{equation}\label{ddimSDE}
\bm{X}_t = \bm{X}_0 + \int_0^t\bm{\mu}(s,\bm{X}_s)\rd s + \int_0^t\bm{\sigma}^2(s,\bm{X}_s)\rd\bm{W}_s,
\end{equation}
$0\leq t<\infty$, where $\bm{\mu}(t,\bm{x})=\{\mu_i(t,x) \}_{1\leq i\leq d}$ is the $(d\times 1)$--drift vector, $\bm{\sigma}(t,\bm{x})=\{\sigma_{ij}(t,x)_{1\leq i\leq d, 1\leq j\leq r}\}$ is the $(d\times r)$- dispersion matrix, $\mu_i(t,x):[0,\infty)\times\mathbb{R}\rightarrow\mathbb{R}$ for $1\leq i\leq d$ and $\sigma_{ij}(t,x):[0,\infty)\times\mathbb{R}\rightarrow\mathbb{R}^+$ for $1\leq i \leq d$, $1\leq j\leq r$ are Borel--measurable functions, and where $W=(\bm{W}_t, \mathscr{F}^W_t; 0\leq t<\infty )$ is the $r$--dimensional Brownian motion in Eq. (\ref{rdimbm}).

One may refer to \cite{ikedawatanabe, karatzasshreve} for the definition of a solution to the SDE \ref{ddimSDE} and the distinction between whether the solution is strong or weak and the type of uniqueness it exhibits.  The following theorem in \cite{karatzasshreve} gives the conditions under which a strong solution with the pathwise uniqueness property exists.
\begin{Theorem}
Suppose the coefficients $\bm{\mu}(t,\bm{x})$ and $\bm{\sigma}(t,\bm{x})$ satisfy the global Lipschitz and linear growth conditions
\begin{eqnarray}
&&\vert\vert\bm{\mu}(t,\bm{x})-\bm{\mu}(t,\bm{y})\vert\vert + \vert\vert\bm{\sigma}(t,\bm{x})- \bm{\sigma}(t,\bm{y})\vert\vert\leq K_1\vert\vert\bm{x}-\bm{y}\vert\vert \\
&&\vert\vert\bm{\mu}(t,\bm{x})\vert\vert^2 + \vert\vert\bm{\sigma}(t,\bm{x})\vert\vert^2\leq K_2\left(1+\vert\vert\bm{x}\vert\vert^2\right)
\end{eqnarray}
for every $0\leq t<\infty$, $\bm{x}\in\mathbb{R}^d$, $\bm{y}\in\mathbb{R}^d$, $0<K_1,K_2<\infty$ and where $\vert\vert\cdot\vert\vert$ denotes the $L^2$ norm.  On some probability space $(\Omega, \mathscr{F},\mathbb{P})$, let $\bm{\xi}$ be an $\mathbb{R}^d$-valued random vector, independent of the $r$--dimensional Brownian motion $W=(\bm{W}_t, \mathscr{F}_t; 0\leq t<\infty )$, and with finite second moment.  Let $\{\mathscr{F}_t\}$ be as in Eq. (\ref{augfiltration}).  Then there exists a continuous, adapted process $X=(\bm{X}_t, \mathscr{F}_t; 0\leq t<\infty)$ which is a strong solution of Eq. (\ref{ddimSDE}) relative to $W$, with initial condition $\bm{\xi}$.  This process is square-integrable.  
\end{Theorem}
\subsection{SDE coefficients of the $h$--transform quantile diffusion}\label{htransformcoefficients}

The dynamics of the $h$--transform quantile diffusion are obtained by Proposition \ref{prop1}, where we can write the density function of the $h$--transform distribution as $$f_{\phi_h}(z) = \expp\left(-0.5(h+1)x_u^2\right)/(\sqrt{2\pi}B(1+h x_u^2))$$ where $u=F_{\phi_h}(z)$ and $\bm{\xi}=(A,B,g)$.  Here, $F_{\phi_h}=Q_{\phi_h}^{-}$ can be computed analytically.  Taking $A=0, B=1$ with no loss of generality for the standardised case, we have the following for the drift function of the $h$--transform quantile diffusion: 

\begin{equation}\begin{split}\label{htransformdriftnottruelaw}
\alpha(t,Z_t) &= \left\{ \mu\left(\cdot,\cdot\right)  f\left(Q\left(F_{\phi_h}(Z_t)\right)\right)+\frac{1}{2}\sigma^2\left(\cdot,\cdot\right) f^{\prime}\left(Q\left(F_{\phi_h}(Z_t)\right)\right) \right\} \left(1+2h\left(\erf^{-}\left(2F_{\phi_h}(Z_t)-1\right)\right)^2\right)
\\[0pt] 
&\qquad\qquad\qquad\times\sqrt{2\pi}\expp\left((h+1)\left(\erf^{-}\left(2F_{\phi_h}(Z_t)-1\right)\right)^2 \right) \\[0pt] 
&\qquad+ \sigma^2\left(\cdot,\cdot\right)f\left(Q\left(F_{\phi_h}(Z_t)\right)\right)^2\left((h+1)\left(1+2h\left(\erf^{-}\left(2F_{\phi_h}(Z_t)-1\right)\right)^2\right)+2h \right)
\\[0pt] 
&\qquad\qquad\qquad\times\sqrt{2}\erf^{-}\left(2F_{\phi_h}(Z_t)-1\right)\pi\expp\left((h+2)\left(\erf^{-}\left(2F_{\phi_h}(Z_t)-1\right)\right)^2 \right)
\end{split}\end{equation}   
for a time--homogeneous $F$ that is the ``false law'' of the driving process, and  
\begin{equation}\begin{split}\label{htransformdrifttruelaw}
\alpha(t,Z_t) &= \sigma^2\left(\cdot,\cdot\right) f^{\prime}\left(t,Q\left(t,F_{\phi_h}(Z_t)\right)\right) \sqrt{2\pi}\left(1+2h\left(\erf^{-}\left(2F_{\phi_h}(Z_t)-1\right)\right)^2\right)
\\
&\qquad\qquad\times\expp\left((h+1)\left(\erf^{-}\left(2F_{\phi_h}(Z_t)-1\right)\right)^2 \right)
 \\[0pt] 
 &\qquad+ \sigma^2\left(\cdot,\cdot\right)f\left(t,Q\left(t,F_{\phi_h}(Z_t)\right)\right)^2 \left((h+1)\left(1+2h\left(\erf^{-}\left(2F_{\phi_h}(Z_t)-1\right)\right)^2\right)+2h \right) 
 \\[0pt] 
 &\qquad\qquad\times\sqrt{2}\erf^{-}\left(2F_{\phi_h}(Z_t)-1\right) \pi\expp\left((h+2)\left(\erf^{-}\left(2F_{\phi_h}(Z_t)-1\right)\right)^2 \right)
\end{split}\end{equation}
when $F$ is the ``true law'' of $(Y_t)$ with $f$ the corresponding transition density.  The volatility function is given by 
\begin{equation}\begin{split}\label{htransformvol}
\widetilde{\sigma}(t,Z_t) &= \sigma\left(\cdot,\cdot\right) f\left(t,Q\left(t,F_{\phi_h}(z)\right)\right)\sqrt{2\pi}\left(1+2h\left(\erf^{-}(2F_{\phi_h}(z)-1)\right)^2\right)\\
&\qquad\qquad\times\expp\left((h+1)\left(\erf^{-}(2F_{\phi_h}(z)-1)\right)^2 \right)
\end{split}\end{equation}

The argument of the drift and diffusion coefficients are given by $\mu(\cdot,\cdot):=\mu(t,Q(t,F_{\phi_h}(z)))$ and $\sigma(\cdot,\cdot):=\sigma(t,Q(t,F_{\phi_h}(z)))$, respectively. In the case where we are using a non time--dependent distribution function in our mapping, we replace $f(t,Q(t,F_{\phi_h}(z)))$ by $f(Q(F_{\phi_h}(z)))$ in Eq. (\ref{htransformvol}).

\subsection{Lipschitz continuity of the $h$--transform quantile diffusion drift and \\ volatility functions}\label{htransformlipconditions}
\begin{Proposition}\label{htransformlipschitzprop}
Let $(Z_t)$ be a $h$--transform quantile process given by Definition \ref{htransformqpdefn}, where the drift coefficient is given by either (\ref{htransformdriftnottruelaw}) or (\ref{htransformdrifttruelaw}), and the volatility coefficient by Eq. (\ref{htransformvol}). Let $(Y_t)_{0\leq t<\infty}$ be a homogeneous driving process satisfying $\rd Y_t = \mu \rd t + \sigma\rd W_t$ for $\mu\in\mathbb{R}$, $\sigma\in\mathbb{R}^+$, and $Y_0=y_0\in\mathbb{R}$.  Then the coefficients of $(Z_t)$ are Lipschitz continuous on $(-\infty, \infty)$ if the density function $f(t,y)$, associated with the law $F_{\bullet}(t,y)$, for all $t\in[t_0,\infty)$ has left and right tail decay to zero, is bounded on its support, and it is such that the following set of conditions is satisfied:
\begin{align}\label{hliplimit1left}
\lim_{x\rightarrow 0^+}\frac{f^{\prime}\left(Q\left(x\right)\right)}{f\left(Q\left(x\right)\right)} &= L_1<\infty, \\[0pt]\label{hliplimit1right}
\lim_{x\rightarrow 1^-}\frac{f^{\prime}\left(Q\left(x\right)\right)}{f\left(Q\left(x\right)\right)} &= L_2<\infty, \\[0pt]\label{hliplimit2left}
\lim_{x\rightarrow 0^+} f\left(Q\left(x\right)\right)(h+1)\erf^{-}\left(2x-1\right)\expp\left(\left(\erf^{-}\left(2x-1\right) \right)^2 \right) &= L_3<\infty, \\[0pt]\label{hliplimit2right}
\lim_{x\rightarrow 1^-} f\left(Q\left(x\right)\right)(h+1)\erf^{-}\left(2x-1\right)\expp\left(\left(\erf^{-}\left(2x-1\right) \right)^2 \right) &= L_4<\infty, \\[0pt]\label{hliplimit3left}
\lim_{x\rightarrow 0^+}\frac{f^{\prime\prime}\left(Q\left(x\right)\right)}{f\left(Q\left(x\right)\right)} &= L_5<\infty, \\[0pt]\label{hliplimit3right}
\lim_{x\rightarrow 1^-}\frac{f^{\prime\prime}\left(Q\left(x\right)\right)}{f\left(Q\left(x\right)\right)} &= L_6<\infty, \\[0pt]\label{hliplimit4left}
\lim_{x\rightarrow 0^+}\frac{f\left(Q\left(x\right)\right)^2\left(\erf^{-}\left(2x-1\right)\right)^6\expp\left(\left(\erf^{-}\left(2x-1\right)\right)^2 \right)}{\left(1+2h\left(\erf^{-}\left(2x-1\right)\right)^2\right)^2} &= L_7 <\infty, \\[0pt]
\label{hliplimit4right}
\lim_{x\rightarrow 1^-}\frac{f\left(Q\left(x\right)\right)^2\left(\erf^{-}\left(2x-1\right)\right)^6\expp\left(\left(\erf^{-}\left(2x-1\right)\right)^2 \right)}{\left(1+2h\left(\erf^{-}\left(2x-1\right)\right)^2\right)^2} &= L_8 <\infty, \\[0pt]\label{hliplimit5left}
\lim_{x\rightarrow 0^+}f\left(Q\left(x\right)\right)^2\left(\erf^{-}\left(2x-1\right)\right)^2\expp\left(2\left(\erf^{-}\left(2x-1\right)\right)^2 \right) &= L_9 <\infty, \\[0pt] \label{hliplimit5right}
\lim_{x\rightarrow 1^-}f\left(Q\left(x\right)\right)^2\left(\erf^{-}\left(2x-1\right)\right)^2\expp\left(2\left(\erf^{-}\left(2x-1\right)\right)^2 \right) &= L_{10} <\infty.
\end{align}
\end{Proposition}

\begin{proof}
The proof is analogous to that of Proposition \ref{gtransformlipschitzprop}. It consists of computing the first partial derivatives of the drift and volatility coefficients given by Eqs (\ref{htransformdriftnottruelaw}) to (\ref{htransformvol}), and finding the conditions that must be satisfied by $f(t,y)$ for all $t\in[t_0,\infty)$ to ensure that these expressions are bounded on the range on which they are differentiable everywhere.
\end{proof}
\end{appendix}

\end{document}